\documentclass{amsart} 

\usepackage{pstricks,color}

\include{pst-plot}

\usepackage{amsmath,amsthm}     
\usepackage{amssymb}            
\usepackage{euscript}           
\usepackage{graphicx,enumerate,calc,lscape,color}
\definecolor{red}{rgb}{1,0,0}
\usepackage[matrix,arrow,curve,frame]{xy}    

\xymatrixcolsep{1.9pc}                          
\xymatrixrowsep{1.9pc}
\newdir{ >}{{}*!/-5pt/\dir{>}}                  

\newtheorem{thm}[subsection]{Theorem}
\newtheorem{defn}[subsection]{Definition}
\newtheorem{prop}[subsection]{Proposition}

\newtheorem{lemma}[subsection]{Lemma}

\theoremstyle{definition}  
\newtheorem{ex}[subsection]{Example}

\newtheorem{remark}[subsection]{Remark}

\newcommand{\tens}              {\otimes}               
\newcommand{\iso}               {\cong}

\newcommand{\field}[1]  {\mathbb #1} 

\newcommand{\bZ}         {\field Z}

\newcommand{\M}         {\field M}

\newcommand{\F}		{\mathbb{F}}
\newcommand{\C}         {\field C}
\newcommand{\Z}         {\field Z}

\DeclareMathOperator{\Spec}{Spec}
\DeclareMathOperator{\Hom}{Hom}

\DeclareMathOperator{\Ext}{Ext}

\DeclareMathOperator{\Gr}{Gr}

\DeclareMathOperator{\Sq}{Sq}

\newcommand{\map}{\rightarrow}

\newcommand{\cl}{\mathrm{cl}}

\newcommand{\tmf}{\textit{tmf}}
\newcommand{\mmf}{\textit{mmf}}

\numberwithin{equation}{subsection}

\newpsobject{taumult}{psline}{linestyle=dashed,dash=0.1 0.2}
\newgray{gridline}{0.9}
\newrgbcolor{multcolor}{0 1 1}
\newpsobject{extn}{psline}{linecolor=red}
\newpsobject{mult}{psline}{linecolor=multcolor}
\newrgbcolor{gsquare}{0.3 0.3 1}
\newrgbcolor{gcube}{1 0 1}

\begin{document}

\title{The cohomology of motivic $A(2)$}

\author{Daniel C.\ Isaksen}

\address{Department of Mathematics\\ Wayne State University\\
Detroit, MI 48202}

\thanks{The 
author was supported by NSF grant DMS0803997.}

\email{isaksen@math.wayne.edu}

\begin{abstract}
Working over an algebraically closed field of characteristic zero,
we compute the cohomology of the subalgebra $A(2)$ of the motivic
Steenrod algebra that is generated by $\Sq^1$, $\Sq^2$, and $\Sq^4$.
The method of calculation is a motivic version of the May spectral
sequence.

Speculatively assuming that there is a ``motivic modular forms" spectrum
with certain properties, we use an Adams-Novikov spectral sequence
to compute the homotopy of such a spectrum at the prime $2$.
\end{abstract}

\maketitle

\section{Introduction}

The purpose of this article is to present some algebraic calculations
that are relevant to motivic homotopy theory.  Recent work \cite{HKO}
\cite{DI} has shown that the Adams and Adams-Novikov spectral sequences
are useful for computing in the motivic stable homotopy category
over an algebraically closed field of characteristic zero, after completion
with respect to the motivic $\F_2$-Eilenberg-Mac Lane spectrum
$H\F_2$ that represents motivic $\F_2$-cohomology.  

This program is built upon 
work of Voevodsky \cite{V1} \cite{V2} \cite{V3}.
Voevodsky has described the motivic $\F_2$-cohomology
of an algebraically closed field of characteristic zero, which we
write as $\M_2$.  He also described the motivic $\F_2$-Steenrod
algebra $A$.  Much like
the classical Steenrod algebra $A_{\cl}$, $A$ is
generated by elements $\Sq^i$, subject to motivic versions of
the Adem relations.

In this article, we take Voevodsky's descriptions as given algebraic
inputs, and we carry out further algebraic computations. 
The motivic Adams spectral sequence takes the cohomology of $A$,
i.e., the ring $\Ext_{A}(\M_2,\M_2)$, as input.
Just as in the classical
situation, it is unlikely that we will ever have a complete description
of this ring.  However, much can be said in low dimensions \cite{DI}.

Classically, one way to approximate difficult calculations
over $A_{\cl}$ is to consider instead the subalgebra $A(2)_{\cl}$ generated
by $\Sq^1$, $\Sq^2$, and $\Sq^4$.  It is possible to give an
explicit but lengthy description of the cohomology of 
$A(2)_{\cl}$ \cite{M2} \cite{IS}.
The main purpose of this article is to carry out the motivic version of
this computation.  Namely, we completely describe the ring
$\Ext_{A(2)}(\M_2,\M_2)$.

The cohomology of
$A(2)_{\cl}$ is the $E_2$-term of an Adams spectral sequence that
converges to the homotopy of $\tmf$ at the prime $2$. 
 The original sources of
this material are mostly unpublished work of Hopkins, Mahowald, and Miller
\cite{H1} \cite{H2}.
One might optimistically hope to construct a motivic spectrum $\mmf$, 
defined over algebraically closed fields of characteristic zero,
whose homotopy is similarly related to the cohomology
of motivic $A(2)$.
This is a potential homotopical application of our algebraic computations.
In this article, we do not discuss any issues related to the existence
of $\mmf$.

Although it is theoretically possible to compute the homotopy of $\tmf$
at the prime $2$
by starting with $\Ext$ groups over $A(2)_{\cl}$ and applying
the Adams spectral sequence, it is difficult in practice.  
The motivic version of this computation is
significantly more difficult.  Even finding the $E_3$-term of
the spectral sequence requires a prohibitive amount of bookkeeping,
and there are a tremendous number of exotic extensions to resolve.

Classically, a more efficient way to compute the homotopy of $\tmf$
at the prime $2$ 
is to start with the cohomology of the elliptic curves Hopf algebroid
and apply the Adams-Novikov spectral sequence.  This computation
is entirely described in \cite{B}, based on work of Hopkins and
Mahowald \cite{HM}.

At the end of this paper, we carry out the motivic version of this
computation.  We make several assumptions about motivic versions
of the elliptic curves Hopf algebroid and the motivic Adams-Novikov
spectral sequence.  Based on these assumptions, we are able to
describe the homotopy of the speculative motivic spectrum $\mmf$
at the prime $2$.

\subsection{Organization of the paper}

Section \ref{sctn:background} gives 
a brief review of the algebraic objects under consideration.
In Section \ref{sctn:May}, we set up a motivic May spectral
sequence that converges to the cohomology of $A(2)$.
We also describe the $E_\infty$-term of this spectral sequence.
In Section \ref{sctn:Ext}, we compute Massey products
and resolve all multiplicative extensions to give a complete
description of the cohomology of $A(2)$ as a ring.
Finally, Section \ref{sctn:mmf} discusses an Adams-Novikov spectral
sequence and describes the homotopy of the speculative motivic
spectrum $\mmf$.


\subsection{Acknowledgments}
We acknowledge the assistance of Robert Bruner.
Section \ref{sctn:May} relies very heavily on unpublished notes of Peter
May on the cohomology of $A(2)_{\cl}$ \cite{M2}.  
Section \ref{sctn:mmf} relies
very heavily on charts of Tilman Bauer \cite{B}.
We also thank Mark
Behrens, Dan Dugger, and Mike Hill for useful conversations.

\section{Background}
\label{sctn:background}

In this section we review the basic algebraic facts about the objects
under consideration.
We are working in categories of bigraded objects.
In a bidegree $(p,q)$, we shall refer to $p$ as the topological degree 
and $q$ as the weight.  The terminology is motivated by the relationship
between motivic cohomology and classical homotopy theory.

\begin{defn}
The bigraded ring $\M_2$ is the polynomial ring $\F_2[\tau]$ on
one generator $\tau$ of bidegree $(0,1)$.
\end{defn}

The relevance of $\M_2$ is that 
it is the motivic $\F_2$-cohomology of an algebraically closed field
of characteristic zero \cite{V1}.

\begin{defn}
The motivic Steenrod algebra $A$ is the $\M_2$-algebra generated by
elements $\Sq^{2k}$ and $\Sq^{2k-1}$ for all $k \geq 1$, of bidegrees
$(2k,k)$ and $(2k-1,k-1)$ respectively, and satisfying the following
relations for $a< 2b$:
\[ 
\Sq^a \Sq^b = 
\sum_{c} \binom{b-1-c}{a-2c} \tau^{?}\Sq^{a+b-c} \Sq^c.
\]
\end{defn}

The expression $\tau^?$ stands for either $1$ or $\tau$.
The distinction is easily determined by consideration of bidegrees.
For example, $\Sq^2 \Sq^2 = \tau \Sq^3 \Sq^1$.

The relevance of $A$ is that it is the ring of motivic $\F_2$-cohomology
operations over an algebraically closed field of characteristic zero
\cite{V2} \cite{V3}.

We consider $\M_2$ as an $A$-module, where
$\Sq^i$ acts trivially on $\M_2$ for $i > 0$ and $\Sq^0$ acts as the identity.
Since $\M_2$ is concentrated in topological degree $0$,
this is the only possible action of $A$ on $\M_2$.

\begin{defn}
The algebra $A(2)$ is the $\M_2$-subalgebra of $A$ generated by 
$\Sq^1$, $\Sq^2$, and $\Sq^4$.  
\end{defn}

\begin{remark}
The algebra $A$ has a Milnor $\M_2$-basis
consisting of elements of the form 
$P^R$, where $R = (r_1, r_2, \ldots)$ ranges over all
finite sequences of non-negative integers.
Just as in the classical case,
$A(2)$ 
has an $\M_2$-basis consisting of elements of the
form $P^R$, where $R = (r_1, r_2, r_3)$, $0 \leq r_1 \leq 7$,
$0 \leq r_2 \leq 3$, and $0 \leq r_1 \leq 1$.
See \cite{DI} for more details on the motivic Milnor basis.
\end{remark}

\subsection{$\Ext$ groups}

We will compute
the tri-graded groups
$\Ext_{A(2)}(\M_2, \M_2)$.

We have
\[
\Ext^{0,(0,*)}_{A(2)}(\M_2,\M_2)=\Hom^{(0,*)}_{A(2)}(\M_2,\M_2)=\M_2.
\]
Here we abuse notation and write $\M_2$ where we really mean the $\M_2$-dual
$\Hom_{\M_2}(\M_2,\M_2)$ of $\M_2$.  
The only important point is that now $\tau$ has bidegree $(0,-1)$.
For fixed $s$ and $t$,
$\Ext_{A(2)}^{s,(t+s,*)}(\M_2,\M_2)$ is a module
over 
$\Ext^{0,(0,*)}_{A(2)}(\M_2,\M_2)$.
In particular, it is an $\M_2$-module and 
therefore decomposes as a sum of
modules of the form $\M_2$ or $\M_2/\tau^k$.  
In Section \ref{subsctn:compare-algebra}, we will explain that 
the free
part coincides with $\Ext$ over the classical version of $A(2)$.

\subsection{Comparison with the classical Steenrod algebra}
\label{subsctn:compare-algebra}

We write $A_{\cl}$ and $A(2)_{\cl}$ for the classical Steenrod
algebra and its subalgebra generated by $\Sq^1$, $\Sq^2$, and $\Sq^4$.

\begin{lemma}
\label{lem:compare-A}
There is a ring isomorphism 
$A[\tau^{-1}] \map A_{\cl} \otimes_{\F_2} \M_2[\tau^{-1}]$
that takes $\Sq^{2k}$ and $\Sq^{2k-1}$ 
to $\tau^{-k} \Sq^{2k}$ and $\tau^{-k} \Sq^{2k-1}$ respectively.
\end{lemma}

In other words, the motivic Steenrod algebra and the classical
Steenrod algebra are essentially the same after inverting $\tau$.
The proof of this lemma is a straightforward algebraic exercise.
Similarly,
there is a ring isomorphism
$A(2)[\tau^{-1}] \map A(2)_{\cl} \otimes_{\F_2} \M_2[\tau^{-1}]$
that takes $\Sq^1$, $\Sq^2$, and $\Sq^4$ to
$\Sq^1$, $\tau^{-1} \Sq^2$, and $\tau^{-2} \Sq^4$ respectively.

\begin{prop} 
\label{prop:compare-ext}
There are isomorphisms of rings
\[\Ext_{A}(\M_2,\M_2)\tens_{\M_2} \M_2[\tau^{-1}] \iso
\Ext_{A_{\cl}}(\F_2,\F_2)\tens_{\F_2} \M_2[\tau^{-1}].
\]
and
\[\Ext_{A(2)}(\M_2,\M_2)\tens_{\M_2}\M_2[\tau^{-1}] \iso
\Ext_{A(2)_{\cl}}(\F_2,\F_2)\tens_{\F_2} \M_2[\tau^{-1}].
\]
\end{prop}

The point of the proposition is that the free $\M_2$-modules in
$\Ext_{A(2)}(\M_2,\M_2)$ correspond precisely to the
free $\F_2$-modules in $\Ext_{A(2)_{\cl}}(\F_2,\F_2)$.
On the other hand, copies of $\M_2/\tau^k$ have no classical
analogue.

The proof of Proposition \ref{prop:compare-ext}
is given in \cite{DI}.  The idea is to use the isomorphism
of Lemma \ref{lem:compare-A} and the flatness of $\M_2 \map \M_2[\tau^{-1}]$.


\section{The motivic May spectral sequence}
\label{sctn:May}

In this section we begin our computation of
$\Ext_{A(2)}(\M_2,\M_2)$ by setting
up a May spectral sequence \cite{M1} and finding the $E_\infty$-term
of this spectral sequence.  
Elements of 
$\Ext_{A(2)}(\M_2,\M_2)$ will be graded in the form
$(s, a, w)$, where $s$ is the stem, i.e., the topological degree minus
the homological degree, $a$ is Adams filtration, i.e., the homological degree,
and $w$ is the weight.

Let $I$ be the two-sided $\M_2$-ideal of
$A(2)$ generated by $\Sq^1$, $\Sq^2$, and $\Sq^4$.
Let $\Gr A(2)$ denote the associated graded algebra $A(2)/I \oplus I/I^2
\oplus I^2/I^3 \oplus\cdots$.  
Elements of $\Gr A(2)$ will be graded in the form $(m, t, w)$,
where $m$ is the May filtration, i.e., the $I$-adic valuation,
$t$ is the topological degree, and $w$ is the weight.

We will first need to compute $\Ext_{\Gr A(2)}(\M_2, \M_2)$.
Elements of this ring will be graded in the form $(m, s, a, w)$,
where $m$ is the May filtration, $s$ is the stem, i.e., the 
topological degree minus the homological degree,
$a$ is the Adams filtration, i.e., the homological degree,
and $w$ is the weight.

\begin{prop}
There is a spectral sequence
\[ 
E_2=\Ext^{(m,s,a,w)}_{\Gr A(2)}(\M_2,M_2) \Rightarrow
\Ext_{A(2)}^{(s,a,w)}(\M_2,\M_2).
\]
\end{prop}

We shall refer to this as the motivic May spectral sequence.
As usual, it can be obtained by filtering the cobar complex by powers
of $I$.  

Let $I_{\cl}$ be the ideal of $A(2)_{\cl}$ 
generated by $\Sq^1$, $\Sq^2$, and $\Sq^4$,
and let $\Gr A(2)_{\cl}$ be the associated graded algebra.

\begin{prop}
\label{prop:assoc-graded}
\mbox{}\par
\begin{enumerate}[(a)]
\item 
The tri-graded algebras $\Gr A(2)$ and
$\Gr A(2)_{\cl} \tens_{\F_2} \F_2[\tau]$ are isomorphic.
\item The quadruply-graded rings
$\Ext_{\Gr A(2)}(\M_2,\M_2)$ and 
$\Ext_{\Gr A(2)_{cl}}(\F_2,\F_2)\tens_{\F_2} \M_2$
are isomorphic.
\end{enumerate}
\end{prop}

The point is that 
the $E_2$-terms of the motivic 
and classical May spectral
sequences are very similar.  Although the motivic and classical
differentials are compatible, they are different in an important
way that we shall describe later.

The proposition can be proved just as in \cite{DI}.  The main point is that
the $\tau$ coefficients in the relations in $A$ appear only on terms
of higher filtration and thus do not affect the associated graded algebra.

The classical ring $\Ext_{\Gr A(2)_{\cl}}(\F_2,\F_2)$ is 
computed in \cite{M2} using the May spectral sequence.  It is also
computed in \cite{IS} using a different method.

By Proposition \ref{prop:assoc-graded} and the results in \cite{M1},
the ring
$\Ext_{\Gr A(2)}(\M_2,\M_2)$ is the cohomology of the differential
graded $\M_2$-algebra
$\M_2[h_{10}, h_{11}, h_{12}, h_{20}, h_{21}, h_{30}]$,
where the differential is described in the following table.
See \cite{DI} for an explanation of the degrees.

\begin{table}[!htbp]
\caption{Differentials for computing $\Ext_{\Gr A(2)}(\F_2,\F_2)$}
\begin{center}
\label{table:Gr-diff}
\begin{tabular}{|l|l|l|}
\hline
$x$ & degree & $d(x)$ \\
\hline
$h_{10}$ & $(1,0,1,0)$ & 0 \\
$h_{11}$ & $(1,1,1,1)$ & 0 \\
$h_{12}$ & $(1,3,1,2)$ & 0 \\
$h_{20}$ & $(2,2,1,1)$ & $h_{10} h_{11}$ \\
$h_{21}$ & $(2,5,1,3)$ & $h_{11} h_{12}$ \\
$h_{30}$ & $(3,6,1,3)$ & $h_{10} h_{21} + h_{20} h_{12}$ \\
\hline
\end{tabular}
\end{center}
\end{table}

It is relatively straightforward to compute
$\Ext_{\Gr A(2)}(\M_2,\M_2)$ from the differentials given above.
The table below lists the generators.

\begin{table}[!htbp]
\caption{Generators for $\Ext_{\Gr A(2)}(\M_2,\M_2)$}
\begin{center}
\label{table:May-E2-gen}
\begin{tabular}{|l|l|l|}
\hline
generator & degree & description in terms of $h_{ij}$ \\
\hline
$h_0$ & $(1,0,1,0)$ & $h_{10}$ \\
$h_1$ & $(1,1,1,1)$ & $h_{11}$ \\
$h_2$ & $(1,3,1,2)$ & $h_{12}$ \\
$b_{20}$ & $(4,4,2,2)$ & $h_{20}^2$ \\
$b_{21}$ & $(4,10,2,6)$ & $h_{21}^2$ \\
$b_{30}$ & $(6,12,2,6)$ & $h_{30}^2$ \\
$h_0(1)$ & $(4,7,2,4)$ & $h_{20} h_{21} + h_{11} h_{30}$ \\
\hline
\end{tabular}
\end{center}
\end{table}

\begin{prop}
\label{prop:May-E2}
The ring 
$\Ext_{\Gr A(2)}(\M_2,\M_2)$ is generated over $\M_2$ by the elements listed in
Table \ref{table:May-E2-gen}, subject to the relations
\begin{enumerate}
\item
$h_0 h_1=0$.
\item
$h_1 h_2=0$.
\item
$h_2 b_{20} = h_0 h_0(1)$.
\item
$h_2 h_0(1) = h_0 b_{21}$.
\item
$h_0(1)^2 = b_{20} b_{21} + h_1^2 b_{30}$.
\end{enumerate}
\end{prop}

The proof of Proposition \ref{prop:May-E2} is a straightforward
lift of the analogous classical computation because of
Proposition \ref{prop:assoc-graded}(b).

\subsection{$E_4$-term of the motivic May spectral sequence}

Having described the $E_2$-term of the motivic May spectral sequence,
i.e., $\Ext_{\Gr A(2)}(\M_2,\M_2)$, we are now ready to compute the 
higher terms.  For dimension reasons, as in the classical case,
the odd differentials must vanish.  In particular, $E_3 = E_4$.

The $d_2$ differentials on the $E_2$-term
are easy to analyze.  They must
be compatible with the $d_2$ differentials in the classical May 
spectral sequence
and they must preserve the weight.  
The following table lists the $d_2$ differentials on all of our generators.
From the data in this table, one can use the Leibniz rule to compute the $d_2$
differential on any element.

\begin{table}[!htbp]
\caption{$d_2$ differentials in the motivic May spectral sequence}
\begin{center}
\label{table:d2}
\begin{tabular}{|l|l||l|l|}
\hline
$x$ & $d_2(x)$ & $x$ & $d_2(x)$ \\
\hline
$h_0$ & 0 &
$b_{20}$ & $\tau h_1^3 + h_0^2 h_2$ \\
$h_1$ & 0 &
$b_{21}$ & $h_2^3$ \\
$h_2$ & 0 &
$b_{30}$ & $\tau h_1 b_{21}$ \\
& & $h_0(1)$ & $h_0 h_2^2$ \\
\hline
\end{tabular}
\end{center}
\end{table}

A straightforward computation now gives the $E_4$-term of the motivic
May spectral sequence.  The generators are listed in the following
table, and the relations are listed in the next theorem.

\begin{table}[!htbp]
\caption{Generators of the $E_4$-term of the motivic May spectral sequence}
\begin{center}
\label{table:May-E4-gen}
\begin{tabular}{|l|l||l|}
\hline
generator & degree & description in $E_2$ \\
\hline
$h_0$ & $(1,0,1,0)$ & $h_0$ \\
$h_1$ & $(1,1,1,1)$ & $h_1$ \\
$h_2$ & $(1,3,1,2)$ & $h_2$  \\
$P$   & $(8,8,4,4)$ & $b_{20}^2$ \\
$c$ & $(5,8,3,5)$ & $h_1 h_0(1)$ \\
$u$   & $(5,11,3,7)$ & $h_1 b_{21}$ \\
$\alpha$ & $(7,12,3,6)$ & $h_0 b_{30}$ \\
$d$ & $(8,14,4,8)$ & $h_1^2 b_{30} + b_{20} b_{21}$ \\
$\nu$ & $(7,15,3,8)$ & $h_2 b_{30}$ \\
$e$ & $(8,17,4,10)$ & $h_0(1) b_{21}$ \\
$g$   & $(8,20,4,12)$ & $b_{21}^2$ \\
$\Delta$ & $(12,24,4,12)$ & $b_{30}^2$ \\
\hline
\end{tabular}
\end{center}
\end{table}

The notation is chosen to be compatible with the standard 
notation for elements in the cohomology of the classical Steenrod algebra.
The elements $c$, $d$, $e$, and $g$ are related to the classical
elements $c_0$, $d_0$, $e_0$, and $g$.
The element $P$ is related to the Adams periodicity operator.
The element $\Delta$ is related to the element in the homotopy 
of $\tmf$ of the same name.
The elements $u$, $\alpha$, and $\nu$ have no analogues and are given
arbitrary names.

\begin{thm}
The $E_4$-term of the motivic May spectral sequence is generated
over $\M_2$ by the elements listed in Table \ref{table:May-E4-gen},
subject to the following relations:
\begin{enumerate}
\item
$h_0 h_1$, 
$h_1 h_2$, 
$h_0^2 h_2 + \tau h_1^3$, 
$h_0 h_2^2$, 
$h_2^3$.
\item
$\tau u$, 
$\tau h_1^2 c$, 
$\tau cd$, 
$\tau ce$,
$\tau cg$.
\item
$h_0^2 \nu + \tau h_1 d$, 
$h_0 h_2 \nu + \tau h_1 e$,
$h_2^2 \nu + \tau h_1 g$.
\item
$h_2 d + h_0 e$,
$h_2 e + h_0 g$,
$h_2 \alpha + h_0 \nu$.
\item
$h_0 c$,
$h_2 c$,
$h_0 u$,
$h_2 u$,
$h_1 \alpha$,
$h_1 \nu$.
\item
$c^2 + h_1^2 d$,
$u^2 + h_1^2 g$,
$c u + h_1^2 e$,
$e^2 + d g$.
\item
$u d + c e$,
$u e + c g$,
$\nu d + \alpha e$,
$\nu e + \alpha g$.
\item
$c \alpha$,
$c \nu$,
$u \alpha$,
$u \nu$.
\item
$\alpha^2 + h_0^2 \Delta$,
$\alpha \nu + h_0 h_2 \Delta$,
$\nu^2 + h_2^2 \Delta$.
\item 
$h_0^2 d + P h_2^2$,
$h_0 \alpha d + P h_2 \nu$,
$d^2 + h_1^4 \Delta + P g$.
\end{enumerate}

\end{thm}

Several observations can be made immediately.  First, the $E_4$-term
contains the polynomial ring $\M_2[P,\Delta]$, and the $E_4$-term
is free as a module over $\M_2[P,\Delta]$.
However, beware that the $E_4$-term is not of the form 
$\M_2[P, \Delta] \otimes_{\M_2} B$, because of the relations in
part (10) of the theorem.

The $E_4$-term also contains the polynomial ring $\M_2[g]$,
but the $E_4$-term is not free as a module over $\M_2[g]$.
For example, $h_0^3 g = 0$, but $h_0^3$ is not zero.
However, the ideal generated by $g$ is free over $\M_2[g]$.

Some parts of the $E_4$-term are depicted in the following chart.
The horizontal axis is the stem, and the vertical axis is the 
Adams filtration.
Solid circles indicate copies of $\M_2$, while open circles
indicate copies of $\M_2/\tau$.
Vertical lines indicate multiplication by $h_0$, lines of slope $1$
indicate multiplication by $h_1$, and lines of slope $\frac{1}{3}$
indicate multiplication by $h_2$.  Dashed lines indicate that
the multiplication hits $\tau$ times a generator.  For example,
the relation $h_0^2 h_2 = \tau h_1^3$ occurs in the 3-stem.

Vertical arrows indicate infinite towers of copies of $\M_2$ connected
by $h_0$ multiplications.
Diagonal arrows indicate infinite towers of copies of $\M_2/\tau$
connected by $h_1$ multiplications.

For legibility, we have omitted most of the multiples of $P$.
We have only shown a few elements in red to express multiplicative
relations with elements that are not multiples of $P$.
We have also omitted the strict multiples of $\Delta$ and $\Delta^2$.

The green parts of the figure consist of elements that are multiples of $g$,
the blue parts consist of elements that are multiples of $g^2$,
and the magenta parts consist of elements that are multiples of $g^3$.
Observe that the blue part is a shifted copy of the green part.
Similarly, if the figure were larger, the magenta part would be another
shifted copy.

\newpage

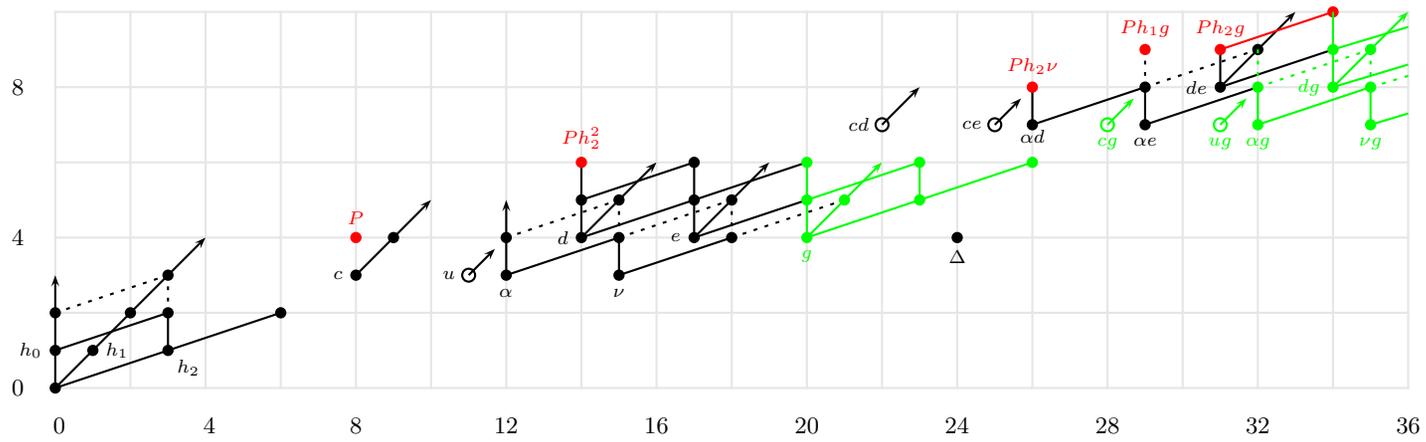
\begin{landscape}

\psset{unit=0.5cm}
\begin{pspicture}(0,0)(36,10)

\psgrid[unit=2,gridcolor=gridline,subgriddiv=0,gridlabelcolor=white](0,0)(18,5)

\small 

\rput(0,-1){ 0}
\rput(4,-1){ 4}
\rput(8,-1){8}
\rput(12,-1){12}
\rput(16,-1){16}
\rput(20,-1){20}
\rput(24,-1){24}
\rput(28,-1){28}
\rput(32,-1){32}
\rput(36,-1){36}

\rput(-1,0){0}
\rput(-1,4){4}
\rput(-1,8){8}


\scriptsize 

\pscircle*(0,0){0.15}
\psline(0,0)(0,1)
\psline(0,0)(1,1)
\psline(0,0)(3,1)
\pscircle*(0,1){0.15}
\psline(0,1)(0,2)
\psline(0,1)(3,2)
\pscircle*(0,2){0.15}
\psline{->}(0,2)(0,3)
\taumult(0,2)(3,3)

\uput[180](0,1){$h_0$}

\pscircle*(1,1){0.15}
\psline(1,1)(2,2)

\uput[0](1,1){$h_1$}

\pscircle*(2,2){0.15}
\psline(2,2)(3,3)

\pscircle*(3,1){0.15}
\psline(3,1)(3,2)
\psline(3,1)(6,2)
\pscircle*(3,2){0.15}
\taumult(3,2)(3,3)
\pscircle*(3,3){0.15}
\psline{->}(3,3)(4,4)

\uput[-45](3,1){$h_2$}

\pscircle*(6,2){0.15}

\pscircle*(8,3){0.15}
\psline(8,3)(9,4)

\uput[180](8,3){$c$}

\pscircle*(9,4){0.15}
\psline{->}(9,4)(10,5)

\pscircle(11,3){0.2}
\psline{->}(11,3)(11.7,3.7)

\uput[180](11,3){$u$}

\pscircle*(12,3){0.15}
\psline(12,3)(12,4)
\psline(12,3)(15,4)
\pscircle*(12,4){0.15}
\psline{->}(12,4)(12,5)
\taumult(12,4)(15,5)

\uput[270](12,3){$\alpha$}

\pscircle*(14,4){0.15}
\psline(14,4)(14,5)
\psline(14,4)(15,5)
\psline(14,4)(17,5)
\pscircle*(14,5){0.15}
\psline(14,5)(14,6)
\psline(14,5)(17,6)

\uput[180](14,4){$d$}

\pscircle*(15,3){0.15}
\psline(15,3)(15,4)
\psline(15,3)(18,4)
\pscircle*(15,4){0.15}
\taumult(15,4)(15,5)
\taumult(15,4)(18,5)
\pscircle*(15,5){0.15}
\psline{->}(15,5)(16,6)

\uput[270](15,3){$\nu$}

\pscircle*(17,4){0.15}
\psline(17,4)(17,5)
\psline(17,4)(18,5)
\psline(17,4)(20,5)
\pscircle*(17,5){0.15}
\psline(17,5)(17,6)
\psline(17,5)(20,6)
\pscircle*(17,6){0.15}

\uput[180](17,4){$e$}

\pscircle*(18,4){0.15}
\taumult(18,4)(18,5)
\taumult(18,4)(21,5)
\pscircle*(18,5){0.15}
\psline{->}(18,5)(19,6)

\pscircle(22,7){0.2}
\psline{->}(22,7)(23,8)

\uput[180](22,7){$cd$}

\pscircle(25,7){0.2}
\psline{->}(25,7)(25.7,7.7)

\uput[180](25,7){$ce$}

\pscircle*(26,7){0.15}
\psline(26,7)(26,8)
\psline(26,7)(29,8)

\uput{0.15}[270](26,7){$\alpha d$}

\pscircle*(29,7){0.15}
\psline(29,7)(29,8)
\psline(29,7)(32,8)
\pscircle*(29,8){0.15}
\taumult(29,8)(29,9)
\taumult(29,8)(32,9)

\uput[270](29,7){$\alpha e$}

\pscircle*(31,8){0.15}
\psline(31,8)(31,9)
\psline(31,8)(32,9)
\psline(31,8)(34,9)

\uput[180](31,8){$de$}

\pscircle*(32,9){0.15}
\psline{->}(32,9)(33,10)

\pscircle*(24,4){0.15}

\uput[270](24,4){$\Delta$}


\psset{linecolor=red}
\color{red}

\pscircle*(8,4){0.15}
\uput[90](8,4){$P$}

\pscircle*(14,6){0.15}
\uput[90](14,6){$P h_2^2$}

\pscircle*(26,8){0.15}
\uput[90](26,8){$P h_2 \nu$}

\pscircle*(29,9){0.15}
\uput[90](29,9){$P h_1 g$}

\pscircle*(31,9){0.15}
\psline(31,9)(34,10)
\uput[90](31,9){$P h_2 g$}

\pscircle*(34,10){0.15}


\psset{linecolor=green}
\color{green}

\pscircle*(20,4){0.15}
\psline(20,4)(20,5)
\psline(20,4)(21,5)
\psline(20,4)(23,5)
\pscircle*(20,5){0.15}
\psline(20,5)(20,6)
\psline(20,5)(23,6)
\pscircle*(20,6){0.15}

\uput[270](20,4){$g$}

\pscircle*(21,5){0.15}
\psline{->}(21,5)(22,6)

\pscircle*(23,5){0.15}
\psline(23,5)(23,6)
\psline(23,5)(26,6)
\pscircle*(23,6){0.15}

\pscircle*(26,6){0.15}

\pscircle(28,7){0.2}
\psline{->}(28,7)(28.7,7.7)

\uput[270](28,7){$cg$}

\pscircle(31,7){0.2}
\psline{->}(31,7)(31.7,7.7)

\uput[270](31,7){$ug$}

\pscircle*(32,7){0.15}
\psline(32,7)(32,8)
\psline(32,7)(35,8)
\pscircle*(32,8){0.15}
\taumult(32,8)(32,9)
\taumult(32,8)(35,9)

\uput[270](32,7){$\alpha g$}

\pscircle*(34,8){0.15}
\psline(34,8)(34,9)
\psline(34,8)(35,9)
\psline(34,8)(36,8.6)
\pscircle*(34,9){0.15}
\psline(34,9)(34,10)
\psline(34,9)(36,9.6)

\uput[180](34,8){$dg$}

\pscircle*(35,7){0.15}
\psline(35,7)(35,8)
\psline(35,7)(36,7.3)
\pscircle*(35,8){0.15}
\taumult(35,8)(35,9)
\taumult(35,8)(36,8.3)
\pscircle*(35,9){0.15}
\psline{->}(35,9)(36,10)

\uput[270](35,7){$\nu g$}

\end{pspicture}

\vskip 2cm

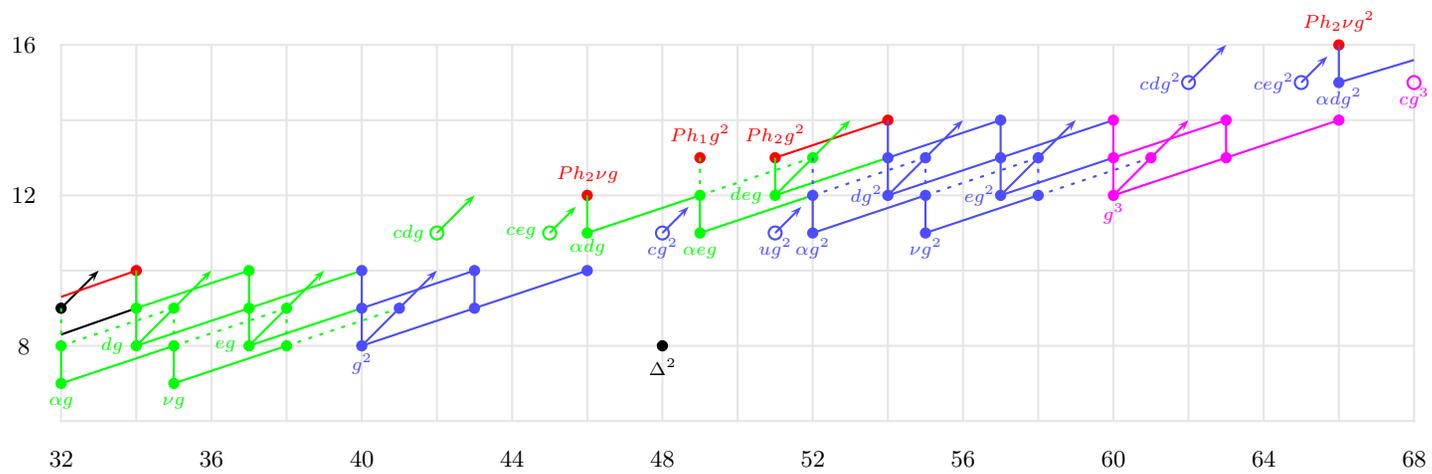
\begin{pspicture}(32,6)(68,16)

\psgrid[unit=2,gridcolor=gridline,subgriddiv=0,gridlabelcolor=white](16,3)(34,8)

\small

\rput(32,5){32}
\rput(36,5){36}
\rput(40,5){40}
\rput(44,5){44}
\rput(48,5){48}
\rput(52,5){52}
\rput(56,5){56}
\rput(60,5){60}
\rput(64,5){64}
\rput(68,5){68}

\rput(31,8){8}
\rput(31,12){12}
\rput(31,16){16}

\scriptsize


\psline(32,8.3)(34,9)

\pscircle*(32,9){0.15}
\psline{->}(32,9)(33,10)

\pscircle*(48,8){0.15}
\uput[270](48,8){$\Delta^2$}


\psset{linecolor=red}
\color{red}

\psline(32,9.3)(34,10)

\pscircle*(34,10){0.15}

\pscircle*(46,12){0.15}
\uput[90](46,12){$Ph_2 \nu g$}

\pscircle*(49,13){0.15}
\uput[90](49,13){$Ph_1 g^2$}

\pscircle*(51,13){0.15}
\psline(51,13)(54,14)
\uput[90](51,13){$Ph_2 g^2$}

\pscircle*(54,14){0.15}

\pscircle*(66,16){0.15}
\uput[90](66,16){$Ph_2 \nu g^2$}


\psset{linecolor=green}
\color{green}

\pscircle*(32,7){0.15}
\psline(32,7)(32,8)
\psline(32,7)(35,8)
\pscircle*(32,8){0.15}
\taumult(32,8)(32,9)
\taumult(32,8)(35,9)

\uput[270](32,7){$\alpha g$}

\pscircle*(34,8){0.15}
\psline(34,8)(34,9)
\psline(34,8)(35,9)
\psline(34,8)(37,9)
\pscircle*(34,9){0.15}
\psline(34,9)(34,10)
\psline(34,9)(37,10)

\uput[180](34,8){$dg$}

\pscircle*(35,7){0.15}
\psline(35,7)(35,8)
\psline(35,7)(38,8)
\pscircle*(35,8){0.15}
\taumult(35,8)(35,9)
\taumult(35,8)(38,9)
\pscircle*(35,9){0.15}
\psline{->}(35,9)(36,10)

\uput[270](35,7){$\nu g$}

\pscircle*(37,8){0.15}
\psline(37,8)(37,9)
\psline(37,8)(38,9)
\psline(37,8)(40,9)
\pscircle*(37,9){0.15}
\psline(37,9)(37,10)
\psline(37,9)(40,10)
\pscircle*(37,10){0.15}

\uput[180](37,8){$eg$}

\pscircle*(38,8){0.15}
\taumult(38,8)(38,9)
\taumult(38,8)(41,9)
\pscircle*(38,9){0.15}
\psline{->}(38,9)(39,10)

\pscircle(42,11){0.2}
\psline{->}(42,11)(43,12)

\uput[180](42,11){$cdg$}

\pscircle(45,11){0.2}
\psline{->}(45,11)(45.7,11.7)

\uput[180](45,11){$ceg$}

\pscircle*(46,11){0.15}
\psline(46,11)(46,12)
\psline(46,11)(49,12)

\uput{0.15}[270](46,11){$\alpha dg$}

\pscircle*(49,11){0.15}
\psline(49,11)(49,12)
\psline(49,11)(52,12)
\pscircle*(49,12){0.15}
\taumult(49,12)(49,13)
\taumult(49,12)(52,13)

\uput[270](49,11){$\alpha eg$}

\pscircle*(51,12){0.15}
\psline(51,12)(51,13)
\psline(51,12)(52,13)
\psline(51,12)(54,13)

\uput[180](51,12){$deg$}

\pscircle*(52,13){0.15}
\psline{->}(52,13)(53,14)


\psset{linecolor=gsquare}
\color{gsquare}

\pscircle*(40,8){0.15}
\psline(40,8)(40,9)
\psline(40,8)(41,9)
\psline(40,8)(43,9)
\pscircle*(40,9){0.15}
\psline(40,9)(40,10)
\psline(40,9)(43,10)
\pscircle*(40,10){0.15}

\uput{0.2}[270](40,8){$g^2$}

\pscircle*(41,9){0.15}
\psline{->}(41,9)(42,10)

\pscircle*(43,9){0.15}
\psline(43,9)(43,10)
\psline(43,9)(46,10)
\pscircle*(43,10){0.15}

\pscircle*(46,10){0.15}

\pscircle(48,11){0.2}
\psline{->}(48,11)(48.7,11.7)

\uput{0.15}[270](48,11){$cg^2$}

\pscircle(51,11){0.2}
\psline{->}(51,11)(51.7,11.7)

\uput{0.15}[270](51,11){$ug^2$}

\pscircle*(52,11){0.15}
\psline(52,11)(52,12)
\psline(52,11)(55,12)
\pscircle*(52,12){0.15}
\taumult(52,12)(52,13)
\taumult(52,12)(55,13)

\uput{0.15}[270](52,11){$\alpha g^2$}

\pscircle*(54,12){0.15}
\psline(54,12)(54,13)
\psline(54,12)(55,13)
\psline(54,12)(57,13)
\pscircle*(54,13){0.15}
\psline(54,13)(54,14)
\psline(54,13)(57,14)

\uput{0.15}[180](54,12){$dg^2$}

\pscircle*(55,11){0.15}
\psline(55,11)(55,12)
\psline(55,11)(58,12)
\pscircle*(55,12){0.15}
\taumult(55,12)(55,13)
\taumult(55,12)(58,13)
\pscircle*(55,13){0.15}
\psline{->}(55,13)(56,14)

\uput{0.15}[270](55,11){$\nu g^2$}

\pscircle*(57,12){0.15}
\psline(57,12)(57,13)
\psline(57,12)(58,13)
\psline(57,12)(60,13)
\pscircle*(57,13){0.15}
\psline(57,13)(57,14)
\psline(57,13)(60,14)
\pscircle*(57,14){0.15}

\uput{0.15}[180](57,12){$eg^2$}

\pscircle*(58,12){0.15}
\taumult(58,12)(58,13)
\taumult(58,12)(61,13)
\pscircle*(58,13){0.15}
\psline{->}(58,13)(59,14)

\pscircle(62,15){0.2}
\psline{->}(62,15)(63,16)

\uput{0.2}[180](62,15){$cdg^2$}

\pscircle(65,15){0.2}
\psline{->}(65,15)(65.7,15.7)

\uput{0.2}[180](65,15){$ceg^2$}

\pscircle*(66,15){0.15}
\psline(66,15)(66,16)
\psline(66,15)(68,15.6)

\uput{0.15}[270](66,15){$\alpha dg^2$}


\psset{linecolor=gcube}

\pscircle*(60,12){0.15}
\psline(60,12)(60,13)
\psline(60,12)(61,13)
\psline(60,12)(63,13)
\pscircle*(60,13){0.15}
\psline(60,13)(60,14)
\psline(60,13)(63,14)
\pscircle*(60,14){0.15}

{\color{gcube} \uput{0.2}[270](60,12){$g^3$}}

\pscircle*(61,13){0.15}
\psline{->}(61,13)(62,14)

\pscircle*(63,13){0.15}
\psline(63,13)(63,14)
\psline(63,13)(66,14)
\pscircle*(63,14){0.15}

\pscircle*(66,14){0.15}

\pscircle(68,15){0.2}

{\color{gcube} \uput{0.15}[270](68,15){$cg^3$}}

\end{pspicture}

\end{landscape}

\subsection{$E_\infty$-term of the motivic May spectral sequence}

Having described the $E_4$-term of the motivic May spectral sequence,
we are now ready to compute the $E_\infty$-term.

As for the $d_2$ differentials, the $d_4$ differentials must be
compatible with the differentials in the classical May spectral sequence.
The only generator that supports a $d_4$ differential is $\Delta$:
\[
d_4(\Delta) = \tau^2 h_2 g.
\]

A straightforward computation now gives the $E_5$-term of the
motivic May spectral sequence.  By inspection, there are no
higher differentials, so $E_5 = E_\infty$.  The generators
of $E_\infty$ are listed in the following table, and the 
relations are listed in the next theorem.

Although $\Delta$ does not survive past $E_4$, we shall write
$\Delta$ for the Massey product operators $\langle \tau^2 g, h_2, - \rangle$
and $\langle h_2, \tau^2 g, - \rangle$.

\begin{table}[!htbp]
\caption{Generators of the $E_\infty$-term of the motivic May spectral sequence}
\begin{center}
\label{table:May-Einfty-gen}
\begin{tabular}{|l|l||l|l|}
\hline
generator & degree & generator & degree \\
\hline
$h_0$ & $(1,0,1,0)$ & $\nu$ & $(7,15,3,8)$  \\
$h_1$ & $(1,1,1,1)$ & $e$ & $(8,17,4,10)$  \\
$h_2$ & $(1,3,1,2)$ & $g$   & $(8,20,4,12)$  \\
$P$   & $(8,8,4,4)$ &
    $\Delta h_1 = \langle \tau^2 g, h_2, h_1 \rangle$ & $(13,25,5,13)$ \\
$c$ & $(5,8,3,5)$ &
    $\Delta c = \langle h_2, \tau^2 g, c\rangle$ & $(17,32,7,17)$ \\
$u$   & $(5,11,3,7)$ &
    $\Delta u = \langle h_2, \tau^2 g, u\rangle$ & $(17,35,7,19)$ \\
$\alpha$ & $(7,12,3,6)$ &
    $\Delta^2$ & $(24,48,8,24)$ \\
$d$ & $(8,14,4,8)$ && \\
\hline
\end{tabular}
\end{center}
\end{table}


\begin{thm}
The $E_\infty$-term of the motivic May spectral sequence is generated
over $\M_2$ by the elements listed in Table \ref{table:May-Einfty-gen},
subject to the relations expressed in the multiplication table
below, as well as the following relations:
\begin{enumerate}
\item
$h_0 h_1$, 
$h_1 h_2$, 
$h_0^2 h_2 + \tau h_1^3$, 
$h_0 h_2^2$, 
$h_2^3$.
\item
$\tau u$, 
$\tau h_1^2 c$, 
$\tau h_1^2 \Delta c$,
$\tau \Delta u$,
$\tau cd$, 
$\tau ce$,
$\tau cg$,
$\tau d \Delta c$,
$\tau e \Delta c$,
$\tau g \Delta c$,
$\tau^2 h_2 g$.
\item
$h_0^2 \nu + \tau h_1 d$, 
$h_0 h_2 \nu + \tau h_1 e$,
$h_2^2 \nu + \tau h_1 g$,
$h_0 \alpha \nu + \tau h_1^2 \Delta h_1$.
\item
$h_2 d + h_0 e$,
$h_2 e + h_0 g$,
$h_2 \alpha + h_0 \nu$.
\item
$h_0 c$,
$h_2 c$,
$h_0 u$,
$h_2 u$,
$h_1 \alpha$,
$h_1 \nu$,
$h_0 \Delta h_1$,
$h_2 \Delta h_1$,
$h_0 \Delta c$,
$h_2 \Delta c$,
$h_0 \Delta u$,
$h_2 \Delta u$,
$h_0 \nu^2$,
$h_2 \nu^2$.
\item
$h_0^2 d + P h_2^2$,
$h_0 \alpha d + P h_2 \nu$,
$\alpha^2 d + P\nu^2$.
\item
$\alpha^2 \nu + \tau d \Delta h_1$,
$\alpha \nu^2 + \tau e \Delta h_1$,
$\nu^3 + \tau g \Delta h_1$,
$\alpha^4 + h_0^4 \Delta^2$.
\end{enumerate}
\end{thm}


\begin{table}[!htbp]
\caption{Multiplication table for $E_\infty$}
\begin{center}
\label{table:Einfty-mult}
\begin{tabular}{|l||l|l|l|l|l|l|l|l|l|}
\hline
& $c$ & $u$ & $\alpha$ & $d$ & $\nu$ & $e$ & $\Delta h_1$ &
    $\Delta c$ & $\Delta u$ \\
\hline
\hline
$c$ & $h_1^2 d$ & $h_1^2 e$ & $0$ & $cd$ & $0$ & $ud$ & $h_1 \Delta c$
    & $h_1 d \Delta h_1$ & $h_1 e \Delta h_1$
\rule{0ex}{2.2ex}\\
\hline
$u$ && $h_1^2 g$ & $0$ & $ce$ & $0$ & $cg$ & $h_1 \Delta u$
    & $h_1 e \Delta h_1$ & $h_1 g \Delta h_1$
\rule{0ex}{2.2ex}\\
\cline{1-1}
\cline{3-10}
$\alpha$ & \multicolumn{2}{c|}{} 
    & $\alpha^2$ & $\alpha d$ & $\alpha \nu$ & $\nu d$ & $0$ & $0$ & $0$
\rule{0ex}{2.2ex}\\
\cline{1-1}
\cline{4-10}
$d$ & \multicolumn{3}{c|}{} 
    & $h_1^3 \Delta h_1 + Pg$ & $\alpha e$ & $de$ & $d \Delta h_1$
    & $d\Delta c$ & $e \Delta c$
\rule{0ex}{2.2ex}\\
\cline{1-1}
\cline{5-10}
$\nu$ & \multicolumn{4}{c|}{} 
    & $\nu^2$ & $\alpha g$ & $0$ & $0$ & $0$
\rule{0ex}{2.2ex}\\
\cline{1-1}
\cline{6-10}
$e$ & \multicolumn{5}{c|}{} 
    & $dg$ & $e\Delta h_1$ & $d \Delta u$ & $g \Delta c$
\rule{0ex}{2.2ex}\\
\cline{1-1}
\cline{7-10}
$\Delta h_1$ & \multicolumn{6}{c|}{} 
    & $h_1^2 \Delta^2$ & $h_1 c \Delta^2$ & $h_1 u \Delta^2$
\rule{0ex}{2.2ex}\\
\cline{1-1}
\cline{8-10}
$\Delta c$ & \multicolumn{7}{c|}{} 
    & $h_1^2 d \Delta^2$ & $h_1^2 e \Delta^2$
\rule{0ex}{2.2ex}\\
\cline{1-1}
\cline{9-10}
$\Delta u$ & \multicolumn{8}{c|}{} 
    & $h_1^2 g \Delta^2$ \rule{0ex}{2.2ex}\\
\hline
\end{tabular}
\end{center}
\end{table}

\newpage

Note that the $E_\infty$-term is free as a module over $\M_2[P,\Delta^2]$.
It is not free as a module over $\M_2[g]$, but the ideal generated
by $g^2$ is free over $\M_2[g]$.

Some parts of the $E_\infty$-term are depicted in the following chart.
The horizontal axis is the stem, and the vertical axis is the Adams
filtration.
Solid circles indicate copies of $\M_2$, open circles indicate
copies of $\M_2/\tau$, open circles with dots indicate copies
of $\M_2/\tau^2$, and open boxes indicate copies of $\M_2/\tau^3$.
Vertical lines indicate multiplication by $h_0$, lines of slope 1 indicate
multiplication by $h_1$, and lines of slope $\frac{1}{3}$ indicate 
multiplication
by $h_2$.  Dashed lines indicate that the multiplication hits $\tau$
times a generator.

Most of the multiples of $P$ are not shown.  A few multiples of $P$ 
are shown in red to express multiplicative relations with elements that
are not multiples of $P$.  Also, multiples of $\Delta^2$ are not shown.

Multiples of $g$ are shown in green, multiples of $g^2$ are shown in blue,
and multiples of $g^3$ are shown in magenta.  If the diagram were larger,
the magenta portion of the diagram would be a shifted copy of the blue portion.

Note the class $de \Delta h_1$ in the 56-stem.  Classically
$de \Delta h_1 = \alpha^3 g$, but motivically we only have
$\alpha^3 g = \tau de \Delta h_1$.  Thus $\tau de \Delta h_1$ is 
a multiple of $g$, but $de \Delta h_1$ is not a multiple of $g$.

\begin{landscape}

\newpsobject{taumult}{psline}{linestyle=dashed,dash=0.1 0.2}

\psset{unit=0.5cm}
\begin{pspicture}(0,0)(36,10)

\psgrid[unit=2,gridcolor=gridline,subgriddiv=0,gridlabelcolor=white](0,0)(18,5)

\small 

\rput(0,-1){ 0}
\rput(4,-1){ 4}
\rput(8,-1){8}
\rput(12,-1){12}
\rput(16,-1){16}
\rput(20,-1){20}
\rput(24,-1){24}
\rput(28,-1){28}
\rput(32,-1){32}
\rput(36,-1){36}

\rput(-1,0){0}
\rput(-1,4){4}
\rput(-1,8){8}


\scriptsize 

\pscircle*(0,0){0.15}
\psline(0,0)(0,1)
\psline(0,0)(1,1)
\psline(0,0)(3,1)
\pscircle*(0,1){0.15}
\psline(0,1)(0,2)
\psline(0,1)(3,2)
\pscircle*(0,2){0.15}
\psline{->}(0,2)(0,3)
\taumult(0,2)(3,3)

\uput[180](0,1){$h_0$}

\pscircle*(1,1){0.15}
\psline(1,1)(2,2)

\uput[0](1,1){$h_1$}

\pscircle*(2,2){0.15}
\psline(2,2)(3,3)

\pscircle*(3,1){0.15}
\psline(3,1)(3,2)
\psline(3,1)(6,2)
\pscircle*(3,2){0.15}
\taumult(3,2)(3,3)
\pscircle*(3,3){0.15}
\psline{->}(3,3)(4,4)

\uput[-45](3,1){$h_2$}

\pscircle*(6,2){0.15}

\pscircle*(8,3){0.15}
\psline(8,3)(9,4)

\uput[180](8,3){$c$}

\pscircle*(9,4){0.15}
\psline{->}(9,4)(10,5)

\pscircle(11,3){0.2}
\psline{->}(11,3)(11.7,3.7)

\uput[180](11,3){$u$}

\pscircle*(12,3){0.15}
\psline(12,3)(12,4)
\psline(12,3)(15,4)
\pscircle*(12,4){0.15}
\psline{->}(12,4)(12,5)
\taumult(12,4)(15,5)

\uput[270](12,3){$\alpha$}

\pscircle*(14,4){0.15}
\psline(14,4)(14,5)
\psline(14,4)(15,5)
\psline(14,4)(17,5)
\pscircle*(14,5){0.15}
\psline(14,5)(14,6)
\psline(14,5)(17,6)

\uput[180](14,4){$d$}

\pscircle*(15,3){0.15}
\psline(15,3)(15,4)
\psline(15,3)(18,4)
\pscircle*(15,4){0.15}
\taumult(15,4)(15,5)
\taumult(15,4)(18,5)
\pscircle*(15,5){0.15}
\psline{->}(15,5)(16,6)

\uput[270](15,3){$\nu$}

\pscircle*(17,4){0.15}
\psline(17,4)(17,5)
\psline(17,4)(18,5)
\psline(17,4)(20,5)
\pscircle*(17,5){0.15}
\psline(17,5)(17,6)
\psline(17,5)(20,6)
\pscircle*(17,6){0.15}

\uput[180](17,4){$e$}

\pscircle*(18,4){0.15}
\taumult(18,4)(18,5)
\taumult(18,4)(21,5)
\pscircle*(18,5){0.15}
\psline{->}(18,5)(19,6)

\pscircle(22,7){0.2}
\psline{->}(22,7)(23,8)

\uput[180](22,7){$cd$}

\pscircle*(24,6){0.15}
\psline{->}(24,6)(24,7)
\taumult(24,6)(27,7)

\uput{0.05}[300](24,6){$\alpha^2$}

\pscircle*(25,5){0.15}
\psline(25,5)(26.2,6)
\pscircle(25,7){0.2}
\psline{->}(25,7)(25.7,7.7)

\uput[270](25,5){$\Delta h_1$}
\uput[180](25,7){$ce$}

\pscircle*(26.2,6){0.15}
\psline(26.2,6)(27,7)
\pscircle*(26,7){0.15}
\psline(26,7)(26,8)
\psline(26,7)(29,8)

\uput{0.15}[200](26,7){$\alpha d$}

\pscircle*(27,6){0.15}
\taumult(27,6)(27,7)
\pscircle*(27,7){0.15}
\psline{->}(27,7)(28,8)

\uput[270](27,6){$\alpha \nu$}

\pscircle*(29,7){0.15}
\psline(29,7)(29,8)
\psline(29,7)(32,8)
\pscircle*(29,8){0.15}
\taumult(29,8)(29,9)
\taumult(29,8)(32,9)

\uput[270](29,7){$\alpha e$}

\pscircle*(30,6){0.15}
\uput{0.2}[270](30,6){$\nu^2$}

\pscircle*(31,8){0.15}
\psline(31,8)(31,9)
\psline(31,8)(32,9)
\psline(31,8)(34,9)

\uput[180](31,8){$de$}

\pscircle*(32,9){0.15}
\psline{->}(32,9)(33,10)

\pscircle*(32.2,7){0.15}
\psline(32.2,7)(33,8)

\uput{0.2}[300](32.2,7){$\Delta c$}

\pscircle*(33,8){0.15}
\psline{->}(33,8)(33.7,8.7)

\pscircle(35.2,7){0.2}
\psline{->}(35.2,7)(36,8)

\uput{0.2}[300](35.2,7){$\Delta u$}

\pscircle*(36,9){0.15}
\psline{->}(36,9)(36,10)

\uput{0.15}[180](36,9){$\alpha^3$}


\psset{linecolor=red}
\color{red}

\pscircle*(8,4){0.15}
\uput[90](8,4){$P$}

\pscircle*(14,6){0.15}
\uput[90](14,6){$P h_2^2$}

\pscircle*(26,8){0.15}
\uput[90](26,8){$P h_2 \nu$}

\pscircle*(29,9){0.15}
\uput[90](29,9){$P h_1 g$}

\pscircle(31,9){0.2}
\pscircle(31,9){0.08}
\psline(31,9)(34,10)
\uput[90](31,9){$P h_2 g$}

\pscircle(34,10){0.2}
\pscircle*(34,10){0.08}


\psset{linecolor=green}
\color{green}

\pscircle*(20,4){0.15}
\psline(20,4)(20,5)
\psline(20,4)(21,5)
\psline(20,4)(23,5)
\pscircle*(20,5){0.15}
\psline(20,5)(20,6)
\psline(20,5)(23,6)
\pscircle*(20,6){0.15}

\uput[270](20,4){$g$}

\pscircle*(21,5){0.15}
\psline{->}(21,5)(22,6)

\pscircle(23,5){0.2}
\pscircle*(23,5){0.08}
\psline(23,5)(23,6)
\psline(23,5)(25.8,6)
\pscircle(23,6){0.2}
\pscircle*(23,6){0.08}

\pscircle(25.8,6){0.2}
\pscircle*(25.8,6){0.08}

\pscircle(28,7){0.2}
\psline{->}(28,7)(28.7,7.7)

\uput[270](28,7){$cg$}

\pscircle(31,7){0.2}
\psline{->}(31,7)(31.7,7.7)

\uput{0.2}[225](31,7){$ug$}

\pscircle*(31.8,7){0.15}
\psline(31.8,7)(32,8)
\psline(31.8,7)(35,8)
\pscircle*(32,8){0.15}
\taumult(32,8)(32,9)
\taumult(32,8)(35,9)

\uput{0.2}[240](31.8,7){$\alpha g$}

\pscircle*(34,8){0.15}
\psline(34,8)(34,9)
\psline(34,8)(35,9)
\psline(34,8)(36,8.6)
\pscircle*(34,9){0.15}
\psline(34,9)(34,10)
\psline(34,9)(36,9.6)

\uput{0.2}[180](34,8){$dg$}

\pscircle*(34.8,7){0.15}
\psline(34.8,7)(35,8)
\psline(34.8,7)(36,7.3)
\pscircle(35,8){0.2}
\pscircle*(35,8){0.08}
\taumult(35,8)(35,9)
\taumult(35,8)(36,8.3)
\rput(35,9){$\Box$}
\psline{->}(35,9)(35.7,9.7)

\uput{0.2}[240](34.8,7){$\nu g$}

\end{pspicture}

\vskip 2cm

\begin{pspicture}(32,6)(68,16)

\psgrid[unit=2,gridcolor=gridline,subgriddiv=0,gridlabelcolor=white](16,3)(34,8)

\small

\rput(32,5){32}
\rput(36,5){36}
\rput(40,5){40}
\rput(44,5){44}
\rput(48,5){48}
\rput(52,5){52}
\rput(56,5){56}
\rput(60,5){60}
\rput(64,5){64}
\rput(68,5){68}

\rput(31,8){8}
\rput(31,12){12}
\rput(31,16){16}

\scriptsize


\psline(32,8.3)(34,9)

\pscircle*(32,9){0.15}
\psline{->}(32,9)(33,10)

\pscircle*(32.2,7){0.15}
\psline(32.2,7)(33,8)

\uput{0.2}[300](32.2,7){$\Delta c$}

\pscircle*(33,8){0.15}
\psline{->}(33,8)(33.7,8.7)

\pscircle(35.2,7){0.2}
\psline{->}(35.2,7)(36,8)

\uput{0.2}[300](35.2,7){$\Delta u$}

\pscircle*(36,9){0.15}
\psline{->}(36,9)(36,10)

\uput{0.15}[180](36,9){$\alpha^3$}

\pscircle*(39,9){0.15}
\psline{->}(39,9)(39.7,9.7)

\pscircle*(41,10){0.15}
\uput{0.05}[225](41,10){$\alpha^2 e$}

\pscircle*(42,9){0.15}
\psline{->}(42,9)(42.7,9.7)
\uput{0.5}[270](42,9){$e \Delta h_1$}

\pscircle(46.2,11){0.2}
\psline{->}(46.2,11)(47,12)
\uput{0.1}[315](46.2,11){$d \Delta c$}

\pscircle(49.2,11){0.2}
\psline{->}(49.2,11)(50,12)
\uput{0.1}[315](49.2,11){$e \Delta c$}

\pscircle*(56,13){0.15}
\psline{->}(56,13)(56.7,13.7)

\pscircle*(48,8){0.15}
\uput[270](48,8){$\Delta^2$}


\psset{linecolor=red}
\color{red}

\psline(32,9.3)(34,10)

\pscircle(34,10){0.2}
\pscircle*(34,10){0.08}

\pscircle(46,12){0.2}
\pscircle*(46,12){0.08}
\uput[90](46,12){$Ph_2 \nu g$}

\rput(49,13){$\Box$}
\uput[90](49,13){$Ph_1 g^2$}

\pscircle(51,13){0.2}
\pscircle*(51,13){0.08}
\psline(51,13)(54,14)
\uput[90](51,13){$Ph_2 g^2$}

\pscircle(54,14){0.2}
\pscircle*(54,14){0.08}

\pscircle(66,16){0.2}
\pscircle*(66,16){0.08}
\uput[90](66,16){$Ph_2 \nu g^2$}


\psset{linecolor=green}
\color{green}

\pscircle*(31.8,7){0.15}
\psline(31.8,7)(32,8)
\psline(31.8,7)(35,8)
\pscircle*(32,8){0.15}
\taumult(32,8)(32,9)
\taumult(32,8)(35,9)

\uput{0.2}[240](31.8,7){$\alpha g$}

\pscircle*(34,8){0.15}
\psline(34,8)(34,9)
\psline(34,8)(35,9)
\psline(34,8)(37,9)
\pscircle*(34,9){0.15}
\psline(34,9)(34,10)
\psline(34,9)(37,10)

\uput{0.2}[180](34,8){$dg$}

\pscircle*(34.8,7){0.15}
\psline(34.8,7)(35,8)
\psline(34.8,7)(38,8)
\pscircle(35,8){0.2}
\pscircle*(35,8){0.08}
\taumult(35,8)(35,9)
\taumult(35,8)(38,9)
\rput(35,9){$\Box$}
\psline{->}(35,9)(35.7,9.7)

\uput{0.2}[240](34.8,7){$\nu g$}

\pscircle*(37,8){0.15}
\psline(37,8)(37,9)
\psline(37,8)(38,9)
\psline(37,8)(40,9)
\pscircle(37,9){0.2}
\pscircle*(37,9){0.08}
\psline(37,9)(37,10)
\psline(37,9)(40,10)
\pscircle(37,10){0.2}
\pscircle*(37,10){0.08}

\uput[180](37,8){$eg$}

\pscircle(38,8){0.2}
\pscircle*(38,8){0.08}
\taumult(38,8)(38,9)
\taumult(38,8)(41,9)
\rput(38,9){$\Box$}
\psline{->}(38,9)(39,10)

\pscircle(42,11){0.2}
\psline{->}(42,11)(43,12)

\uput[180](42,11){$cdg$}

\pscircle*(44,10){0.15}
\uput{0.05}[225](44,10){$\alpha^2 g$}

\pscircle*(45,9){0.15}
\psline{->}(45,9)(45.7,9.7)
\pscircle(45,11){0.2}
\psline{->}(45,11)(45.7,11.7)

\uput{0.2}[270](45,9){$g \Delta h_1$}
\uput[180](45,11){$ceg$}

\pscircle*(45.8,11){0.15}
\psline(45.8,11)(46,12)
\psline(45.8,11)(49,12)

\uput{0.15}[270](45.8,11){$\alpha dg$}

\pscircle*(47,10){0.15}
\uput{0.2}[270](47,10){$\alpha \nu g$}

\pscircle*(48.8,11){0.15}
\psline(48.8,11)(49,12)
\psline(48.8,11)(52,12)
\pscircle(49,12){0.2}
\pscircle*(49,12){0.08}
\taumult(49,12)(49,13)
\taumult(49,12)(52,13)

\uput[270](48.8,11){$\alpha eg$}

\pscircle*(50,10){0.15}
\uput{0.2}[270](50,10){$\nu^2 g$}

\pscircle*(51,12){0.15}
\psline(51,12)(51,13)
\psline(51,12)(52,13)
\psline(51,12)(54,13)

\uput{0.15}[180](51,12){$deg$}

\pscircle(52.2,11){0.2}
\psline{->}(52.2,11)(53,12)
\uput{0.2}[315](52.2,11){$g\Delta c$}

\pscircle(55.2,11){0.2}
\psline{->}(55.2,11)(56,12)
\uput{0.2}[315](55.2,11){$g\Delta u$}

\rput(52,13){$\Box$}
\psline{->}(52,13)(53,14)

\pscircle*(59,13){0.15}
\psline{->}(59,13)(59.7,13.7)

\pscircle*(61,14){0.15}
\uput{0.0}[240](61,14){$\alpha^2 eg$}

\pscircle*(62,13){0.15}
\psline{->}(62,13)(62.7,13.7)
\uput{0.5}[270](62,13){$eg \Delta h_1$}

\pscircle(66.2,15){0.2}
\psline{->}(66.2,15)(67,16)
\uput{0.1}[315](66.2,15){$dg \Delta c$}


\psset{linecolor=gsquare}
\color{gsquare}

\pscircle*(40,8){0.15}
\psline(40,8)(40,9)
\psline(40,8)(41,9)
\psline(40,8)(43,9)
\pscircle(40,9){0.2}
\pscircle*(40,9){0.08}
\psline(40,9)(40,10)
\psline(40,9)(43,10)
\pscircle(40,10){0.2}
\pscircle*(40,10){0.08}

\uput{0.2}[270](40,8){$g^2$}

\rput(41,9){$\Box$}
\psline{->}(41,9)(42,10)

\pscircle(43,9){0.2}
\pscircle*(43,9){0.08}
\psline(43,9)(43,10)
\psline(43,9)(46,10)
\pscircle(43,10){0.2}
\pscircle*(43,10){0.08}

\pscircle(46,10){0.2}
\pscircle*(46,10){0.08}

\pscircle(48,11){0.2}
\psline{->}(48,11)(48.7,11.7)

\uput{0.15}[270](48,11){$cg^2$}

\pscircle(51,11){0.2}
\psline{->}(51,11)(51.7,11.7)

\uput{0.15}[270](51,11){$ug^2$}

\pscircle*(51.8,11){0.15}
\psline(51.8,11)(52,12)
\psline(51.8,11)(55,12)
\pscircle(52,12){0.2}
\pscircle*(52,12){0.08}
\taumult(52,12)(52,13)
\taumult(52,12)(55,13)

\uput{0.25}[270](51.8,11){$\alpha g^2$}

\pscircle*(54,12){0.15}
\psline(54,12)(54,13)
\psline(54,12)(55,13)
\psline(54,12)(57,13)
\pscircle(54,13){0.2}
\pscircle*(54,13){0.08}
\psline(54,13)(54,14)
\psline(54,13)(57,14)

\uput{0.15}[180](54,12){$dg^2$}

\pscircle*(54.8,11){0.15}
\psline(54.8,11)(55,12)
\psline(54.8,11)(58,12)
\pscircle(55,12){0.2}
\pscircle*(55,12){0.08}
\taumult(55,12)(55,13)
\taumult(55,12)(58,13)
\rput(55,13){$\Box$}
\psline{->}(55,13)(56,14)

\uput{0.2}[240](54.8,11){$\nu g^2$}

\pscircle*(57,12){0.15}
\psline(57,12)(57,13)
\psline(57,12)(58,13)
\psline(57,12)(60,13)
\pscircle(57,13){0.2}
\pscircle*(57,13){0.08}
\psline(57,13)(57,14)
\psline(57,13)(60,14)
\pscircle(57,14){0.2}
\pscircle*(57,14){0.08}

\uput{0.15}[180](57,12){$eg^2$}

\pscircle(58,12){0.2}
\pscircle*(58,12){0.08}
\taumult(58,12)(58,13)
\taumult(58,12)(61,13)
\rput(58,13){$\Box$}
\psline{->}(58,13)(59,14)

\pscircle(62,15){0.2}
\psline{->}(62,15)(63,16)

\uput{0.2}[180](62,15){$cdg^2$}

\pscircle*(64,14){0.15}
\uput{0.0}[255](64,14){$\alpha^2 g^2$}

\pscircle*(65,13){0.15}
\psline{->}(65,13)(65.7,13.7)
\pscircle(65,15){0.2}
\psline{->}(65,15)(65.7,15.7)

\uput{0.2}[270](65,13){$g^2 \Delta h_1$}
\uput{0.2}[180](65,15){$ceg^2$}

\pscircle*(65.8,15){0.15}
\psline(65.8,15)(66,16)
\psline(65.8,15)(68,15.6)

\uput{0.2}[240](65.8,15){$\alpha dg^2$}

\pscircle*(67,14){0.15}
\uput{0.2}[270](67,14){$\alpha \nu g^2$}


\newrgbcolor{gcube}{1 0 1}
\psset{linecolor=gcube}

\pscircle*(60,12){0.15}
\psline(60,12)(60,13)
\psline(60,12)(61,13)
\psline(60,12)(63,13)
\pscircle(60,13){0.2}
\pscircle*(60,13){0.08}
\psline(60,13)(60,14)
\psline(60,13)(63,14)
\pscircle(60,14){0.2}
\pscircle*(60,14){0.08}

{\color{gcube} \uput{0.2}[270](60,12){$g^3$}}
\rput(60,14){}

{\color{gcube} \rput(61,13){$\Box$}}
\psline{->}(61,13)(62,14)

\pscircle(63,13){0.2}
\pscircle*(63,13){0.08}
\psline(63,13)(63,14)
\psline(63,13)(66,14)
\pscircle(63,14){0.2}
\pscircle*(63,14){0.08}

\pscircle(66,14){0.2}
\pscircle*(66,14){0.08}

\pscircle(68,15){0.2}

{\color{gcube} \uput{0.15}[270](68,15){$cg^3$}}

\rput(39,9){}

{\color{black}
\uput{0.2}[270](39,9){$d \Delta h_1$} 
\uput{0.2}[270](56,13){$de\Delta h_1$} }

\rput(59,13){}

{\color{green}
\uput{0.2}[270](59,13){$dg \Delta h_1$} }

\end{pspicture}

\end{landscape}

\section{Computation of $\Ext_{A(2)}(\M_2,\M_2)$}
\label{sctn:Ext}

With the $E_\infty$-term of the motivic May spectral sequence in hand
from the previous section, we have nearly computed
$\Ext_{A(2)}(\M_2,\M_2)$.  It only remains to resolve some hidden
extensions in the multiplicative structure.

\subsection{Massey products}

In this section we compute some Massey products 
in $\Ext_{A(2)}(\M_2,\M_2)$
that will be
needed to resolve the extension problems.
All Massey products that we consider have zero indeterminacy.

\begin{lemma}
In $\Ext_{A(2)}(\M_2,\M_2)$,
\begin{enumerate}
\item
$\tau h_1^2 = \langle h_0, h_1, h_0 \rangle$.
\item
$h_0 h_2 = \langle h_1, h_0, h_1 \rangle$.
\item
$h_2^2 = \langle h_1, h_2, h_1 \rangle$.
\end{enumerate}
\end{lemma}

There are several ways of understanding these formulas.  The simplest is
to note that these formulas are already known in 
$\Ext_{A(2)_{\cl}}(\F_2,\F_2)$, and comparison
gives the motivic formulas also.
Another approach
is to compute a minimal resolution for $A(2)$ in low dimensions and 
to compute Yoneda products and Massey products
explicitly in terms of chain maps.

\begin{lemma}
In $\Ext_{A(2)}(\M_2,\M_2)$,
\begin{enumerate}
\item
$\tau d = \langle h_0, h_1, \alpha \rangle$.
\item
$\tau e = \langle h_0, h_1, \nu \rangle$.
\item
$\tau e = \langle h_2, h_1, \alpha \rangle$.
\item
$\tau g = \langle h_2, h_1, \nu \rangle$.
\end{enumerate}
\end{lemma}

\begin{proof}
For the first formula, we have
\[
h_1 \langle h_0, h_1, \alpha \rangle = 
\langle h_1, h_0, h_1 \rangle \alpha = h_0 h_2 \alpha = \tau h_1 d.
\]
In the relevant dimension, multiplication by $h_1$ is injective,
so $\tau d = \langle h_0, h_1, \alpha \rangle$.

The proofs of the other formulas are similar.
\end{proof}

\begin{lemma}
In $\Ext_{A(2)}(\M_2,\M_2)$,
\begin{enumerate}
\item
$c = \langle h_1, h_0, h_2^2 \rangle$.
\item
$u = \langle h_1, h_2, h_2^2 \rangle$.
\item
$\alpha = \langle h_0, h_1, h_2, \tau h_2^2 \rangle$.
\item
$h_0 d = \langle \tau h_1 c, h_1, h_2 \rangle$.
\item
$\alpha e = \langle \tau g, c, h_0 \rangle$.
\item
$\alpha g = \langle \tau g, u, h_0 \rangle$.
\end{enumerate}
\end{lemma}

\begin{proof}
The formulas are already true in the $E_4$-term of the motivic
May spectral sequence.  They can be verified by computations with
explicit cocycles.  There are no extension problems to resolve in passing
from $E_4$ to $\Ext_{A(2)}(\M_2, \M_2)$.
\end{proof}

We remind the reader that $\Delta h_1$, $\Delta c$, and $\Delta u$
are defined as follows.

\begin{defn}
In $\Ext_{A(2)}(\M_2, \M_2)$,
\begin{enumerate}
\item
$\Delta h_1 = \langle \tau^2 g, h_2, h_1 \rangle$.
\item
$\Delta c = \langle h_2, \tau^2 g, c \rangle$.
\item
$\Delta u = \langle h_2, \tau^2 g, u \rangle$.
\end{enumerate}
\end{defn}

\subsection{Multiplicative Extensions}

\begin{lemma}
In $\Ext_{A(2)}(\M_2,\M_2)$, the elements
$\tau \Delta u$, $\tau d \Delta c$, $\tau e \Delta c$,
and $\tau g \Delta c$ are all zero.
\end{lemma}

\begin{proof}
For the first element, note that
$\tau \langle u, \tau^2 g, h_2 \rangle \subseteq
\langle 0, \tau^2 g, h_2 \rangle = 0$.

For the second element,
\[
\langle h_2, \tau^2 g, c \rangle \tau d =
h_2 \langle \tau^2 g, c, \tau d \rangle,
\]
which equals zero since $\langle \tau^2 g, c, \tau d \rangle$
is zero for dimension reasons.
The argument for the last two elements is similar.
\end{proof}

\begin{lemma}
In $\Ext_{A(2)}(\M_2,\M_2)$, we have 
\begin{enumerate}
\item
$d \Delta u = e \Delta c$.
\item
$e \Delta u = g \Delta c$.
\end{enumerate}
\end{lemma}

\begin{proof}
We know from $E_\infty$ that $e \Delta u$ equals either
$g \Delta c$ or $g \Delta c + \tau \alpha g^2$.
We have already shown that $\tau \cdot g \Delta c$ 
and $\tau \cdot e \Delta u$ are both zero,
but $\tau^3 \alpha g^2$ is non-zero.  Therefore,
$e \Delta u = g \Delta c$.

The proof of the second formula is similar.
\end{proof}

\begin{lemma}
In $\Ext_{A(2)}(\M_2,\M_2)$, we have the following formulas:
\begin{enumerate}
\item
$h_0 \Delta c = \tau h_0 \alpha g$.
\item
$h_0 \Delta u = \tau h_0 \nu g$.
\item
$h_2 \Delta c = \tau h_0 \nu g$.
\item
$h_2 \Delta u = \tau h_2 \nu g$.
\end{enumerate}
\end{lemma}

\begin{proof}
For the first formula,
\[
h_0 \langle c, \tau^2 g, h_2 \rangle =
\langle h_0, c, \tau g \rangle \tau h_2 =
\tau h_2 \alpha e = \tau h_0 \alpha g.
\]
The proof of the second formula is similar.

The third and fourth formulas follow easily by multiplying
the first and second formulas by $h_2$, using that multiplication
by $h_0$ is injective in the relevant dimensions.
\end{proof}

\begin{lemma}
In $\Ext_{A(2)}(\M_2,\M_2)$, we have the following formulas:
\begin{enumerate}
\item
$\alpha c = \tau h_0^2 g$.
\item
$\nu c = \tau h_0 h_2 g$.
\item
$\alpha u = \tau h_0 h_2 g$.
\item
$\nu u = \tau h_2^2 g$.
\item
$\alpha \Delta h_1 = \tau^3 eg$.
\item
$\nu \Delta h_1 = \tau^3 g^2$.
\end{enumerate}
\end{lemma}

\begin{proof}
For the first formula,
\[
\alpha \langle h_1, h_0, h_2^2 \rangle = 
\langle \alpha, h_1, h_0 \rangle h_2^2 = \tau d h_2^2 = \tau h_0^2 g.
\]
The other calculations are similar.
\end{proof}

\begin{lemma}
In $\Ext_{A(2)}(\M_2,\M_2)$, the products
$\alpha \Delta c$,
$\alpha \Delta u$,
$\nu \Delta c$, and
$\nu \Delta u$ are all zero.
\end{lemma}

\begin{proof}
For the first product,
\[
\alpha \Delta c = 
\langle h_0, h_1, h_2, \tau h_2^2 \rangle \Delta c \subseteq
\langle h_0, h_1, h_2, \tau^3 h_1 e g \rangle = 
\langle h_0, h_1, h_2, 0 \rangle = 0.
\]
The second product vanishes for similar reasons.

For the third product,
\[
\tau \nu \Delta c =
\langle h_2, \tau^2 g, c \rangle \tau \nu =
h_2 \langle \tau^2 g, c, \tau \nu \rangle,
\]
which is zero because 
$\langle \tau^2 g, c, \tau \nu \rangle$ is zero for dimension reasons.
Therefore, $\nu \Delta c$ does not equal $\tau \alpha \nu g$
since $\tau^2 \alpha \nu g$ is not zero.
The fourth product vanishes for similar reasons.
\end{proof}

\begin{lemma}
In $\Ext_{A(2)}(\M_2,\M_2)$, we have the following formulas:
\begin{enumerate}
\item
$(\Delta h_1)^2 = h_1^2 \Delta^2 + \tau^2 \nu^2 g$.
\item
$\Delta h_1 \cdot \Delta c = h_1 c \Delta^2$.
\item
$\Delta h_1 \cdot \Delta u = h_1 u \Delta^2$.
\item
$(\Delta c)^2 = h_1^2 d \Delta^2$.
\item
$\Delta c \cdot \Delta u = h_1^2 e \Delta^2$.
\item
$(\Delta u)^2 = h_1^2 g \Delta^2$.
\end{enumerate}
\end{lemma}

\begin{proof}
For the first formula, 
we know from $E_\infty$ that 
$(\Delta h_1)^2$ equals either $h_1^2 \Delta^2$ or
$h_1^2 \Delta^2 + \tau^2 \nu^2 g$.
We have already shown that $\nu (\Delta h_1)^2 = \tau^3 g^2 \Delta h_1$
is non-zero, but $\nu \cdot h_1^2 \Delta^2$ is zero.
It follows that 
$(\Delta h_1)^2 = h_1^2 \Delta^2 + \nu^2 g$.

For the second formula, we know from $E_\infty$ that
$\Delta h_1 \cdot \Delta c$ equals either $h_1 c \Delta^2$ or
$h_1 c \Delta^2 + \tau^4 eg^2$.  
We have already shown that 
$\tau g \cdot \Delta h_1 \Delta c$ and 
$\tau g \cdot h_1 c \Delta^2$ are both zero, but
$\tau g \cdot \tau^4 eg^3$ is non-zero.

For the third formula, we know from $E_\infty$ that 
$\Delta h_1 \cdot \Delta u$ equals $h_1 u \Delta^2$
or $h_1 u \Delta^2 + \tau^4 g^3$.
We have already shown that
$\tau \cdot \Delta h_1 \Delta u$ and
$\tau \cdot h_1 u \Delta^2$ are zero,
but $\tau \cdot \tau^4 g^3$ is non-zero.

For the fourth formula, we know from $E_\infty$ that 
$(\Delta c)^2$ equals $h_1^2 d \Delta^2$ or
$h_1^2 d \Delta^2 + \tau^2 \alpha^2 g^2$.
Observe that 
\[
\tau \cdot (\Delta c)^2 =
\langle h_2, \tau^2 g, c \rangle \tau \Delta c =
h_2 \langle \tau^2 g, c, \tau \Delta c \rangle,
\]
which must be zero because 
$\langle \tau^2 g, c, \tau \Delta c \rangle$ is zero for
dimension reasons.
Also, $\tau \cdot h_1^2 d \Delta^2$ is zero, but
$\tau \cdot \tau^2 \alpha^2 g^2$ is non-zero.

For the fifth formula, we know from $E_\infty$ that 
$\Delta c \cdot \Delta u$ equals $h_1^2 e \Delta^2$
or $h_1^2 e \Delta^2 + \tau^2 \alpha \nu g^2$.
We already know that $\tau \cdot (\Delta c \cdot \Delta u)$
and $\tau \cdot h_1^2 e \Delta^2$ are zero, but
$\tau \cdot \tau^2 \alpha \nu g$ is non-zero.

For the sixth formula, we know from $E_\infty$ that 
$(\Delta u)^2$ equals $h_1^2 g \Delta^2$
or $h_1^2 g \Delta^2 + \tau^2 \nu^2 g^2$.
We already know that $\tau \cdot (\Delta u)^2$ and
$\tau \cdot h_1^2 g \Delta^2$ are zero,
but $\tau \cdot \tau^2 \nu^2 g^2$ is non-zero.
\end{proof}

\begin{lemma}
In $\Ext_{A(2)}(\M_2,\M_2)$, we have 
$\tau h_1^2 \Delta c = \tau^2 h_0 d g$ but
$h_1^2 \Delta c$ does not equal $\tau h_0 dg$.
\end{lemma}

\begin{proof}
From $E_\infty$, we know that $h_1 \cdot h_1^2 \Delta c$ is non-zero,
but $h_1 \cdot \tau h_0 dg$ is zero.  Therefore,
$h_1^2 \Delta c$ and $\tau h_0 dg$ cannot be equal.

From $E_\infty$, we also know that $h_1 \Delta c$
equals $c \Delta h_1$.  Therefore,
\[
\tau h_1^2 \Delta c = \tau h_1 c \Delta h_1 =
\tau h_1 c \langle h_1, h_2, \tau^2 g \rangle =
\langle \tau h_1 c, h_1, h_2 \rangle \tau^2 g.
\]
It follows that $\tau h_1^2 \Delta c = \tau^2 h_0 dg$.
\end{proof}

\begin{lemma}
In $\Ext_{A(2)}(\M_2,\M_2)$, we have 
$\alpha^4 = h_0^4 \Delta^2 + \tau^4 P g^2$.
\end{lemma}

\begin{proof}
From $E_\infty$, we know that $\alpha^4$ equals either 
$h_0^4 \Delta^2$ or $h_0^4 \Delta^2 + \tau^4 P g^2$.
We have already shown that $\Delta h_1 \cdot \alpha^4$ equals
\[
\alpha^3 \cdot \tau^3 eg = 
\tau^3 \alpha^2 e \cdot \nu e =
\tau^4 d \Delta h_1 \cdot dg =
\tau^4 g \Delta h_1 ( Pg + h_1^3 \Delta h_1) =
\tau^4 P g^2 \Delta h_1.
\]
However, $\Delta h_1 \cdot h_0^4 \Delta^2$ equals zero.
It follows that $\alpha^4$ must equal 
$h_0^4 \Delta^2 + \tau^4 P g^2$.
\end{proof}

\subsection{Ring structure of $\Ext_{A(2)}(\M_2,\M_2)$}

Finally, we assemble the calculations made in the previous sections
to compute $\Ext_{A(2)}(\M_2,\M_2)$ as a ring.  The generators are listed
in the following table.  All possible multiplicative extensions have
been resolved in the previous section.

\begin{table}[!htbp]
\caption{Generators of $\Ext_{A(2)}$}
\begin{center}
\label{table:Ext-gen}
\begin{tabular}{|l|l|l|l|}
\hline
generator & degree & generator & degree\\
\hline
$h_0$ & $(0,1,0)$ &
    $\nu$ & $(15,3,8)$  \\
$h_1$ & $(1,1,1)$ &
    $e$ & $(17,4,10)$  \\
$h_2$ & $(3,1,2)$ &
    $g$   & $(20,4,12)$  \\
$P$   & $(8,4,4)$ &
    $\Delta h_1$ & $(25,5,13)$ \\
$c$ & $(8,3,5)$   &
    $\Delta c$ & $(32,7,17)$ \\
$u$   & $(11,3,7)$ &
    $\Delta u$ & $(35,7,19)$ \\
$\alpha$ & $(12,3,6)$ &
    $\Delta^2$ & $(48,8,24)$ \\
$d$ & $(14,4,8)$  && \\
\hline
\end{tabular}
\end{center}
\end{table}

\begin{thm}
The ring $\Ext_{A(2)}(\M_2,\M_2)$ is generated
over $\M_2$ by the elements listed in Table \ref{table:Ext-gen},
subject to the relations expressed in the multiplication table
below, as well as the following relations:
\begin{enumerate}
\item
$h_0 h_1$, 
$h_1 h_2$, 
$h_0^2 h_2 + \tau h_1^3$, 
$h_0 h_2^2$, 
$h_2^3$.
\item
$\tau u$, 
$\tau h_1^2 c$, 
$\tau h_1^2 \Delta c + \tau^2 h_0 dg$,
$\tau \Delta u$,
$\tau cd$, 
$\tau ce$,
$\tau cg$,
$\tau d \Delta c$,
$\tau e \Delta c$,
$\tau g \Delta c$,
$\tau^2 h_2 g$.
\item
$h_0^2 \nu + \tau h_1 d$, 
$h_0 h_2 \nu + \tau h_1 e$,
$h_2^2 \nu + \tau h_1 g$,
$h_0 \alpha \nu + \tau h_1^2 \Delta h_1$.
\item
$h_2 d + h_0 e$,
$h_2 e + h_0 g$,
$h_2 \alpha + h_0 \nu$.
\item
$h_0 c$,
$h_2 c$,
$h_0 u$,
$h_2 u$,
$h_1 \alpha$,
$h_1 \nu$,
$h_0 \Delta h_1$,
$h_2 \Delta h_1$,
$h_0 \Delta c + \tau h_0 \alpha g$,
$h_2 \Delta c + \tau h_0 \nu g$,
$h_0 \Delta u + \tau h_0 \nu g$,
$h_2 \Delta u + \tau h_2 \nu g$,
$h_0 \nu^2$,
$h_2 \nu^2$.
\item
$h_0^2 d + P h_2^2$,
$h_0 \alpha d + P h_2 \nu$,
$\alpha^2 d + P\nu^2$.
\item
$\alpha^2 \nu + \tau d \Delta h_1$,
$\alpha \nu^2 + \tau e \Delta h_1$,
$\nu^3 + \tau g \Delta h_1$,
$\alpha^4 + h_0^4 \Delta^2 + \tau^4 Pg^2$.
\end{enumerate}
\end{thm}

\newpage


\begin{table}[!htbp]
\caption{Multiplication table for $\Ext_{A(2)}$}
\begin{center}
\label{table:Ext-mult}
\begin{tabular}{|l||l|l|l|l|l|}
\hline
& $c$ & $u$ & $\alpha$ & $d$ & $\nu$ \\
\hline
\hline
$c$ & $h_1^2 d$ & $h_1^2 e$ & $\tau h_0^2 g$ & $cd$ & $\tau h_0 h_2 g$
\rule{0ex}{2.2ex}\\
\hline
$u$ && $h_1^2 g$ & $\tau h_0 h_2 g$ & $ce$ & $\tau h_2^2 g$ 
\rule{0ex}{2.2ex}\\
\cline{1-1}
\cline{3-6}
$\alpha$ & \multicolumn{2}{c|}{} 
    & $\alpha^2$ & $\alpha d$ & $\alpha \nu$ 
\rule{0ex}{2.2ex}\\
\cline{1-1}
\cline{4-6}
$d$ & \multicolumn{3}{c|}{} 
    & $h_1^3 \Delta h_1 + Pg$ & $\alpha e$ 
\rule{0ex}{2.2ex}\\
\cline{1-1}
\cline{5-6}
$\nu$ & \multicolumn{4}{c|}{} 
    & $\nu^2$ 
\rule{0ex}{2.2ex}\\
\hline
\multicolumn{6}{c}{} \\
\cline{1-5}
& $e$ & $\Delta h_1$ &
    $\Delta c$ & $\Delta u$ \\
\cline{1-5}
\cline{1-5}
$c$ & $ud$ & $h_1 \Delta c$ & $h_1 d \Delta h_1$ & $h_1 e \Delta h_1$
\rule{0ex}{2.2ex}\\
\cline{1-5}
$u$ & $cg$ & $h_1 \Delta u$ & $h_1 e \Delta h_1$ & $h_1 g \Delta h_1$
\rule{0ex}{2.2ex}\\
\cline{1-5}
$\alpha$ & $\nu d$ & $\tau^3 eg$ & $0$ & $0$
\rule{0ex}{2.2ex}\\
\cline{1-5}
$d$ & $de$ & $d \Delta h_1$ & $d\Delta c$ & $e \Delta c$
\rule{0ex}{2.2ex}\\
\cline{1-5}
$\nu$ & $\alpha g$ & $\tau^3 g^2$ & $0$ & $0$
\rule{0ex}{2.2ex}\\
\cline{1-5}
$e$ & $dg$ & $e\Delta h_1$ & $d \Delta u$ & $g \Delta c$
\rule{0ex}{2.2ex}\\
\cline{1-5}
$\Delta h_1$ & \multicolumn{1}{c|}{} 
    & $h_1^2 \Delta^2 + \tau^2 \nu^2 g$ & $h_1 c \Delta^2$ & $h_1 u \Delta^2$
\rule{0ex}{2.2ex}\\
\cline{1-1}
\cline{3-5}
$\Delta c$ & \multicolumn{2}{c|}{} 
    & $h_1^2 d \Delta^2$ & $h_1^2 e \Delta^2$
\rule{0ex}{2.2ex}\\
\cline{1-1}
\cline{4-5}
$\Delta u$ & \multicolumn{3}{c|}{} 
    & $h_1^2 g \Delta^2$ \rule{0ex}{2.2ex}\\
\cline{1-5}
\end{tabular}
\end{center}
\end{table}

Some parts of $\Ext_{A(2)}(\M_2,\M_2)$ are depicted in the following
chart.  
The notation is the same as in the chart for the $E_\infty$-term
of the motivic May spectral sequence.

The most interesting difference occurs with $\Delta c$.  Note
that $\Delta c$ supports exotic multiplications by $h_0$ and $h_2$.
Note also that $\tau h_1^2 \Delta c$ is no longer zero; this 
is an exotic $\tau$-extension.  In fact, $h_1^2 \Delta c$
is the sum $(h_1^2 \Delta c + \tau h_0 dg) + \tau h_0 dg$,
where the first term is killed by $\tau$ 
and the second is killed by no power of $\tau$.

Similarly to $\Delta c$, the classes $\Delta u$, $d \Delta c$, and $e \Delta c$
all support exotic multiplications by $h_0$ and $h_2$.  

\subsection{An Adams spectral sequence?}

At this point, it is natural to wonder whether the cohomology of $A(2)$
is the $E_2$-term of an Adams spectral sequence that converges to the
homotopy of some motivic spectrum that is analogous to $\tmf$.

Assuming that such a spectral sequence exists, it is possible to 
determine the $d_2$-differentials by comparison to calculations over 
the full motivic Steenrod algebra $A$ and by comparison to calculations
over the classical $A(2)_{\cl}$.

For example, comparison to the full motivic Steenrod algebra implies
that $d_2(e) = h_1^2 d$ \cite{DI}.  Then the relation $ud = ce$ implies
that $d_2(u) = h_1^2 c$.  Similar kinds of arguments allow one to
compute the entire $d_2$ differential.  

We have not pursued this computation any further because there is a 
simpler way of making the same speculative calculation that we present
in the following section.

\begin{landscape}

\newpsobject{taumult}{psline}{linestyle=dashed,dash=0.1 0.2}

\psset{unit=0.5cm}
\begin{pspicture}(0,0)(36,10)

\psgrid[unit=2,gridcolor=gridline,subgriddiv=0,gridlabelcolor=white](0,0)(18,5)

\small 

\rput(0,-1){ 0}
\rput(4,-1){ 4}
\rput(8,-1){8}
\rput(12,-1){12}
\rput(16,-1){16}
\rput(20,-1){20}
\rput(24,-1){24}
\rput(28,-1){28}
\rput(32,-1){32}
\rput(36,-1){36}

\rput(-1,0){0}
\rput(-1,4){4}
\rput(-1,8){8}


\scriptsize 

\pscircle*(0,0){0.15}
\psline(0,0)(0,1)
\psline(0,0)(1,1)
\psline(0,0)(3,1)
\pscircle*(0,1){0.15}
\psline(0,1)(0,2)
\psline(0,1)(3,2)
\pscircle*(0,2){0.15}
\psline{->}(0,2)(0,3)
\taumult(0,2)(3,3)

\uput[180](0,1){$h_0$}

\pscircle*(1,1){0.15}
\psline(1,1)(2,2)

\uput[0](1,1){$h_1$}

\pscircle*(2,2){0.15}
\psline(2,2)(3,3)

\pscircle*(3,1){0.15}
\psline(3,1)(3,2)
\psline(3,1)(6,2)
\pscircle*(3,2){0.15}
\taumult(3,2)(3,3)
\pscircle*(3,3){0.15}
\psline{->}(3,3)(4,4)

\uput[-45](3,1){$h_2$}

\pscircle*(6,2){0.15}

\pscircle*(8,3){0.15}
\psline(8,3)(9,4)

\uput[180](8,3){$c$}

\pscircle*(9,4){0.15}
\psline{->}(9,4)(10,5)

\pscircle(11,3){0.2}
\psline{->}(11,3)(11.7,3.7)

\uput[180](11,3){$u$}

\pscircle*(12,3){0.15}
\psline(12,3)(12,4)
\psline(12,3)(15,4)
\pscircle*(12,4){0.15}
\psline{->}(12,4)(12,5)
\taumult(12,4)(15,5)

\uput[270](12,3){$\alpha$}

\pscircle*(14,4){0.15}
\psline(14,4)(14,5)
\psline(14,4)(15,5)
\psline(14,4)(17,5)
\pscircle*(14,5){0.15}
\psline(14,5)(14,6)
\psline(14,5)(17,6)

\uput[180](14,4){$d$}

\pscircle*(15,3){0.15}
\psline(15,3)(15,4)
\psline(15,3)(18,4)
\pscircle*(15,4){0.15}
\taumult(15,4)(15,5)
\taumult(15,4)(18,5)
\pscircle*(15,5){0.15}
\psline{->}(15,5)(16,6)

\uput[270](15,3){$\nu$}

\pscircle*(17,4){0.15}
\psline(17,4)(17,5)
\psline(17,4)(18,5)
\psline(17,4)(20,5)
\pscircle*(17,5){0.15}
\psline(17,5)(17,6)
\psline(17,5)(20,6)
\pscircle*(17,6){0.15}

\uput[180](17,4){$e$}

\pscircle*(18,4){0.15}
\taumult(18,4)(18,5)
\taumult(18,4)(21,5)
\pscircle*(18,5){0.15}
\psline{->}(18,5)(19,6)

\pscircle(22,7){0.2}
\psline{->}(22,7)(23,8)

\uput[180](22,7){$cd$}

\pscircle*(24,6){0.15}
\psline{->}(24,6)(24,7)
\taumult(24,6)(27,7)

\uput{0.05}[300](24,6){$\alpha^2$}

\pscircle*(25,5){0.15}
\psline(25,5)(26.2,6)
\pscircle(25,7){0.2}
\psline{->}(25,7)(25.7,7.7)

\uput[270](25,5){$\Delta h_1$}
\uput[180](25,7){$ce$}

\pscircle*(26.2,6){0.15}
\psline(26.2,6)(27,7)
\pscircle*(26,7){0.15}
\psline(26,7)(26,8)
\psline(26,7)(29,8)

\uput{0.15}[200](26,7){$\alpha d$}

\pscircle*(27,6){0.15}
\taumult(27,6)(27,7)
\pscircle*(27,7){0.15}
\psline{->}(27,7)(28,8)

\uput[270](27,6){$\alpha \nu$}

\pscircle*(29,7){0.15}
\psline(29,7)(29,8)
\psline(29,7)(32,8)
\pscircle*(29,8){0.15}
\taumult(29,8)(29,9)
\taumult(29,8)(32,9)

\uput[270](29,7){$\alpha e$}

\pscircle*(30,6){0.15}
\uput{0.2}[270](30,6){$\nu^2$}

\pscircle*(31,8){0.15}
\psline(31,8)(31,9)
\psline(31,8)(32,9)
\psline(31,8)(33.8,9)

\uput[180](31,8){$de$}

\pscircle*(32,9){0.15}
\psline{->}(32,9)(33,10)

\pscircle*(32.3,7){0.15}
\psline(32.3,7)(33,8)
\taumult(32.3,7)(32,8)
\taumult(32.3,7)(35,8)

\uput{0.2}[300](32.3,7){$\Delta c$}

\pscircle*(33,8){0.15}
\psline(33,8)(34.2,9)
\taumult(33,8)(33.8,9)

\pscircle(34.2,9){0.2}
\psline{->}(34.2,9)(35,10)

\pscircle(35.3,7){0.2}
\psline{->}(35.3,7)(36,8)
\taumult(35.3,7)(35,8)
\taumult(35.3,7)(36,7.3)

\uput{0.2}[300](35.3,7){$\Delta u$}

\pscircle*(36,9){0.15}
\psline{->}(36,9)(36,10)

\uput{0.15}[180](36,9){$\alpha^3$}


\psset{linecolor=red}
\color{red}

\pscircle*(8,4){0.15}
\uput[90](8,4){$P$}

\pscircle*(14,6){0.15}
\uput[90](14,6){$P h_2^2$}

\pscircle*(26,8){0.15}
\uput[90](26,8){$P h_2 \nu$}

\pscircle*(29,9){0.15}
\uput[90](29,9){$P h_1 g$}

\pscircle(31,9){0.2}
\pscircle(31,9){0.08}
\psline(31,9)(34,10)
\uput[90](31,9){$P h_2 g$}

\pscircle(34,10){0.2}
\pscircle*(34,10){0.08}


\psset{linecolor=green}
\color{green}

\pscircle*(20,4){0.15}
\psline(20,4)(20,5)
\psline(20,4)(21,5)
\psline(20,4)(23,5)
\pscircle*(20,5){0.15}
\psline(20,5)(20,6)
\psline(20,5)(23,6)
\pscircle*(20,6){0.15}

\uput[270](20,4){$g$}

\pscircle*(21,5){0.15}
\psline{->}(21,5)(22,6)

\pscircle(23,5){0.2}
\pscircle*(23,5){0.08}
\psline(23,5)(23,6)
\psline(23,5)(25.8,6)
\pscircle(23,6){0.2}
\pscircle*(23,6){0.08}

\pscircle(25.8,6){0.2}
\pscircle*(25.8,6){0.08}

\pscircle(28,7){0.2}
\psline{->}(28,7)(28.7,7.7)

\uput[270](28,7){$cg$}

\pscircle(31,7){0.2}
\psline{->}(31,7)(31.7,7.7)

\uput{0.2}[225](31,7){$ug$}

\pscircle*(31.7,7){0.15}
\psline(31.7,7)(32,8)
\psline(31.7,7)(35,8)
\pscircle*(32,8){0.15}
\taumult(32,8)(32,9)
\taumult(32,8)(35,9)

\uput{0.2}[270](31.7,7){$\alpha g$}

\pscircle*(34,8){0.15}
\psline(34,8)(33.8,9)
\psline(34,8)(35,9)
\psline(34,8)(36,8.6)
\pscircle*(33.8,9){0.15}
\psline(33.8,9)(34,10)
\psline(33.8,9)(36,9.6)

\uput{0.2}[180](34,8){$dg$}

\pscircle*(34.7,7){0.15}
\psline(34.7,7)(35,8)
\psline(34.7,7)(36,7.3)
\pscircle(35,8){0.2}
\pscircle*(35,8){0.08}
\taumult(35,8)(35,9)
\taumult(35,8)(36,8.3)
\rput(35,9){$\Box$}
\psline{->}(35,9)(35.7,9.7)

\uput{0.2}[240](34.7,7){$\nu g$}

\end{pspicture}

\vskip 2cm

\begin{pspicture}(32,6)(68,16)

\psgrid[unit=2,gridcolor=gridline,subgriddiv=0,gridlabelcolor=white](16,3)(34,8)

\small

\rput(32,5){32}
\rput(36,5){36}
\rput(40,5){40}
\rput(44,5){44}
\rput(48,5){48}
\rput(52,5){52}
\rput(56,5){56}
\rput(60,5){60}
\rput(64,5){64}
\rput(68,5){68}

\rput(31,8){8}
\rput(31,12){12}
\rput(31,16){16}

\scriptsize


\psline(32,8.3)(33.8,9)

\pscircle*(32,9){0.15}
\psline{->}(32,9)(33,10)

\pscircle*(32.3,7){0.15}
\psline(32.3,7)(33,8)
\taumult(32.3,7)(32,8)
\taumult(32.3,7)(35,8)

\uput{0.2}[300](32.3,7){$\Delta c$}

\pscircle*(33,8){0.15}
\psline(33,8)(34.2,9)
\taumult(33,8)(33.8,9)

\pscircle(34.2,9){0.2}
\psline{->}(34.2,9)(35,10)

\pscircle(35.3,7){0.2}
\psline{->}(35.3,7)(36,8)
\taumult(35.3,7)(35,8)
\taumult(35.3,7)(38,8)

\uput{0.2}[300](35.3,7){$\Delta u$}

\pscircle*(36,9){0.15}
\psline{->}(36,9)(36,10)

\uput{0.15}[180](36,9){$\alpha^3$}

\pscircle*(39,9){0.15}
\psline{->}(39,9)(39.7,9.7)

\pscircle*(41,10){0.15}
\uput{0.05}[225](41,10){$\alpha^2 e$}

\pscircle*(42,9){0.15}
\psline{->}(42,9)(42.7,9.7)
\uput{0.5}[270](42,9){$e \Delta h_1$}

\pscircle(46.3,11){0.2}
\psline{->}(46.3,11)(47,12)
\taumult(46.3,11)(46,12)
\taumult(46.3,11)(49,12)
\uput{0.1}[315](46.3,11){$d \Delta c$}

\pscircle(49.3,11){0.2}
\psline{->}(49.3,11)(50,12)
\taumult(49.3,11)(49,12)
\taumult(49.3,11)(52,12)
\uput{0.1}[315](49.3,11){$e \Delta c$}

\pscircle*(56,13){0.15}
\psline{->}(56,13)(56.7,13.7)

\pscircle*(48,8){0.15}
\uput[270](48,8){$\Delta^2$}


\psset{linecolor=red}
\color{red}

\psline(32,9.3)(34,10)

\pscircle(34,10){0.2}
\pscircle*(34,10){0.08}

\pscircle(46,12){0.2}
\pscircle*(46,12){0.08}
\uput[90](46,12){$Ph_2 \nu g$}

\rput(49,13){$\Box$}
\uput[90](49,13){$Ph_1 g^2$}

\pscircle(51,13){0.2}
\pscircle*(51,13){0.08}
\psline(51,13)(54,14)
\uput[90](51,13){$Ph_2 g^2$}

\pscircle(54,14){0.2}
\pscircle*(54,14){0.08}

\pscircle(66,16){0.2}
\pscircle*(66,16){0.08}
\uput[90](66,16){$Ph_2 \nu g^2$}


\psset{linecolor=green}
\color{green}

\pscircle*(31.7,7){0.15}
\psline(31.7,7)(32,8)
\psline(31.7,7)(35,8)
\pscircle*(32,8){0.15}
\taumult(32,8)(32,9)
\taumult(32,8)(35,9)

\uput{0.2}[240](31.7,7){$\alpha g$}

\pscircle*(34,8){0.15}
\psline(34,8)(33.8,9)
\psline(34,8)(35,9)
\psline(34,8)(37,9)
\pscircle*(33.8,9){0.15}
\psline(33.8,9)(34,10)
\psline(33.8,9)(37,10)

\uput{0.2}[180](34,8){$dg$}

\pscircle*(34.7,7){0.15}
\psline(34.7,7)(35,8)
\psline(34.7,7)(38,8)
\pscircle(35,8){0.2}
\pscircle*(35,8){0.08}
\taumult(35,8)(35,9)
\taumult(35,8)(38,9)
\rput(35,9){$\Box$}
\psline{->}(35,9)(35.7,9.7)

\uput{0.2}[240](34.7,7){$\nu g$}

\pscircle*(37,8){0.15}
\psline(37,8)(37,9)
\psline(37,8)(38,9)
\psline(37,8)(40,9)
\pscircle(37,9){0.2}
\pscircle*(37,9){0.08}
\psline(37,9)(37,10)
\psline(37,9)(40,10)
\pscircle(37,10){0.2}
\pscircle*(37,10){0.08}

\uput[180](37,8){$eg$}

\pscircle(38,8){0.2}
\pscircle*(38,8){0.08}
\taumult(38,8)(38,9)
\taumult(38,8)(41,9)
\rput(38,9){$\Box$}
\psline{->}(38,9)(39,10)

\pscircle(42,11){0.2}
\psline{->}(42,11)(43,12)

\uput[180](42,11){$cdg$}

\pscircle*(44,10){0.15}
\uput{0.05}[225](44,10){$\alpha^2 g$}

\pscircle*(45,9){0.15}
\psline{->}(45,9)(45.7,9.7)
\pscircle(45,11){0.2}
\psline{->}(45,11)(45.7,11.7)

\uput{0.2}[270](45,9){$g \Delta h_1$}
\uput[180](45,11){$ceg$}

\pscircle*(45.7,11){0.15}
\psline(45.7,11)(46,12)
\psline(45.7,11)(49,12)

\uput{0.15}[270](45.7,11){$\alpha dg$}

\pscircle*(47,10){0.15}
\uput{0.2}[270](47,10){$\alpha \nu g$}

\pscircle*(48.7,11){0.15}
\psline(48.7,11)(49,12)
\psline(48.7,11)(52,12)
\pscircle(49,12){0.2}
\pscircle*(49,12){0.08}
\taumult(49,12)(49,13)
\taumult(49,12)(52,13)

\uput[285](48.7,11){$\alpha eg$}

\pscircle*(50,10){0.15}
\uput{0.2}[270](50,10){$\nu^2 g$}

\pscircle*(51,12){0.15}
\psline(51,12)(51,13)
\psline(51,12)(52,13)
\psline(51,12)(54,13)

\uput{0.15}[180](51,12){$deg$}

\pscircle(52.3,11){0.2}
\psline{->}(52.3,11)(53,12)
\taumult(52.3,11)(52,12)
\taumult(52.3,11)(55,12)
\uput{0.2}[315](52.3,11){$g\Delta c$}

\pscircle(55.3,11){0.2}
\psline{->}(55.3,11)(56,12)
\taumult(55.3,11)(55,12)
\taumult(55.3,11)(58,12)
\uput{0.2}[315](55.3,11){$g\Delta u$}

\rput(52,13){$\Box$}
\psline{->}(52,13)(53,14)

\pscircle*(59,13){0.15}
\psline{->}(59,13)(59.7,13.7)

\pscircle*(61,14){0.15}
\uput{0.0}[240](61,14){$\alpha^2 eg$}

\pscircle*(62,13){0.15}
\psline{->}(62,13)(62.7,13.7)
\uput{0.5}[270](62,13){$eg \Delta h_1$}

\pscircle(66.3,15){0.2}
\psline{->}(66.3,15)(67,16)
\taumult(66.3,15)(66,16)
\taumult(66.3,15)(68,15.6)
\uput{0.2}[300](66.3,15){$dg \Delta c$}


\psset{linecolor=gsquare}
\color{gsquare}

\pscircle*(40,8){0.15}
\psline(40,8)(40,9)
\psline(40,8)(41,9)
\psline(40,8)(43,9)
\pscircle(40,9){0.2}
\pscircle*(40,9){0.08}
\psline(40,9)(40,10)
\psline(40,9)(43,10)
\pscircle(40,10){0.2}
\pscircle*(40,10){0.08}

\uput{0.2}[270](40,8){$g^2$}

\rput(41,9){$\Box$}
\psline{->}(41,9)(42,10)

\pscircle(43,9){0.2}
\pscircle*(43,9){0.08}
\psline(43,9)(43,10)
\psline(43,9)(46,10)
\pscircle(43,10){0.2}
\pscircle*(43,10){0.08}

\pscircle(46,10){0.2}
\pscircle*(46,10){0.08}

\pscircle(48,11){0.2}
\psline{->}(48,11)(48.7,11.7)

\uput{0.15}[270](48,11){$cg^2$}

\pscircle(51,11){0.2}
\psline{->}(51,11)(51.7,11.7)

\uput{0.15}[270](51,11){$ug^2$}

\pscircle*(51.7,11){0.15}
\psline(51.7,11)(52,12)
\psline(51.7,11)(55,12)
\pscircle(52,12){0.2}
\pscircle*(52,12){0.08}
\taumult(52,12)(52,13)
\taumult(52,12)(55,13)

\uput{0.25}[270](51.7,11){$\alpha g^2$}

\pscircle*(54,12){0.15}
\psline(54,12)(54,13)
\psline(54,12)(55,13)
\psline(54,12)(57,13)
\pscircle(54,13){0.2}
\pscircle*(54,13){0.08}
\psline(54,13)(54,14)
\psline(54,13)(57,14)

\uput{0.15}[180](54,12){$dg^2$}

\pscircle*(54.7,11){0.15}
\psline(54.7,11)(55,12)
\psline(54.7,11)(58,12)
\pscircle(55,12){0.2}
\pscircle*(55,12){0.08}
\taumult(55,12)(55,13)
\taumult(55,12)(58,13)
\rput(55,13){$\Box$}
\psline{->}(55,13)(56,14)

\uput{0.2}[240](54.7,11){$\nu g^2$}

\pscircle*(57,12){0.15}
\psline(57,12)(57,13)
\psline(57,12)(58,13)
\psline(57,12)(60,13)
\pscircle(57,13){0.2}
\pscircle*(57,13){0.08}
\psline(57,13)(57,14)
\psline(57,13)(60,14)
\pscircle(57,14){0.2}
\pscircle*(57,14){0.08}

\uput{0.15}[180](57,12){$eg^2$}

\pscircle(58,12){0.2}
\pscircle*(58,12){0.08}
\taumult(58,12)(58,13)
\taumult(58,12)(61,13)
\rput(58,13){$\Box$}
\psline{->}(58,13)(59,14)

\pscircle(62,15){0.2}
\psline{->}(62,15)(63,16)

\uput{0.2}[180](62,15){$cdg^2$}

\pscircle*(64,14){0.15}
\uput{0.0}[255](64,14){$\alpha^2 g^2$}

\pscircle*(65,13){0.15}
\psline{->}(65,13)(65.7,13.7)
\pscircle(65,15){0.2}
\psline{->}(65,15)(65.7,15.7)

\uput{0.2}[270](65,13){$g^2 \Delta h_1$}
\uput{0.2}[180](65,15){$ceg^2$}

\pscircle*(65.7,15){0.15}
\psline(65.7,15)(66,16)
\psline(65.7,15)(68,15.6)

\uput{0.2}[240](65.7,15){$\alpha dg^2$}

\pscircle*(67,14){0.15}
\uput{0.2}[270](67,14){$\alpha \nu g^2$}


\newrgbcolor{gcube}{1 0 1}
\psset{linecolor=gcube}

\pscircle*(60,12){0.15}
\psline(60,12)(60,13)
\psline(60,12)(61,13)
\psline(60,12)(63,13)
\pscircle(60,13){0.2}
\pscircle*(60,13){0.08}
\psline(60,13)(60,14)
\psline(60,13)(63,14)
\pscircle(60,14){0.2}
\pscircle*(60,14){0.08}

{\color{gcube} \uput{0.2}[270](60,12){$g^3$}}

\rput(60,14){}

{\color{gcube} \rput(61,13){$\Box$}}
\psline{->}(61,13)(62,14)

\pscircle(63,13){0.2}
\pscircle*(63,13){0.08}
\psline(63,13)(63,14)
\psline(63,13)(66,14)
\pscircle(63,14){0.2}
\pscircle*(63,14){0.08}

\pscircle(66,14){0.2}
\pscircle*(66,14){0.08}

\pscircle(68,15){0.2}

{\color{gcube} \uput{0.15}[270](68,15){$cg^3$}}


\rput(39,9){}

{\color{black}
\uput{0.2}[270](39,9){$d \Delta h_1$}
\uput{0.2}[270](56,13){$de\Delta h_1$}
}

\rput(59,13){}

{\color{green}
\uput{0.2}[270](59,13){$dg \Delta h_1$}
}

\end{pspicture}
\end{landscape}

\section{Adams-Novikov spectral sequence for ``motivic modular forms''}
\label{sctn:mmf}

In this section, we make a leap of faith and assume that there exists
a motivic spectrum $\mmf$, called ``motivic modular forms", defined over
$\Spec \C$.  We assume that the topological realization of this 
motivic spectrum is the classical spectrum $\tmf$.
We also assume that the homotopy of $\mmf$
can be computed by an Adams-Novikov spectral sequence whose $E_2$-term
is the cohomology of a version of the elliptic curves Hopf algebroid.

Recall that the elliptic curves Hopf algebroid localized at the prime $2$
has the form 
$(\tilde{A}_{\cl},\tilde{\Gamma}_{\cl})$, 
where $\tilde{A}_{\cl}$ is the ring $\Z_{(2)}[a_1, a_3, a_4, a_6]$
and $\tilde{\Gamma}_{\cl} = \tilde{A}_{\cl}[s,t]$ \cite{B}.

We write $\M_{(2)}$ for the ring $\Z_{(2)}[\tau]$, i.e., the 
$2$-local motivic cohomology of a point.

\begin{defn}
The motivic elliptic curves Hopf algebroid localized at $2$ is
$(\tilde{A}, \tilde{\Gamma})$, where
$\tilde{A} = \M_{(2)} [a_1, a_3, a_4, a_6]$ and
$\tilde{\Gamma} = \tilde{A}[s,t]$.  The bidegree of $a_i$ is $(2i,i)$,
while the bidegrees of $s$ and $t$ are $(2,1)$ and $(6,3)$ respectively.
\end{defn}

The structure maps of $(\tilde{A}, \tilde{\Gamma})$ are defined to be
compatible with the classical structure maps.

One may compute the cohomology of $(\tilde{A},\tilde{\Gamma})$
as in \cite[Section 7]{B}.  It turns out that the weights introduce
no new complications.  The answer is essentially the same 
as the classical answer, as shown in the chart on page 30 of \cite{B}.
In fact, the motivic answer is equal to the classical answer tensored
over $\Z_{(2)}$ with $\M_{(2)}$.
Also, the generators are assigned weights, but these are easy to 
determine.  At location $(s,t)$ of the chart,
the weight of the generator is $\frac{s+t}{2}$.

We now have the $E_2$-term of a (speculative) spectral sequence
for computing the homotopy groups of $\mmf$ at the prime $2$.  Next
we consider differentials and compute the $E_\infty$-term.

We assume that the topological realization
of $\mmf$ is $\tmf$.  This implies that the motivic differentials
are algebraically compatible with the classical differentials,
i.e., they are equal after inverting $\tau$.
Using this fact, 
it turns out that the motivic differentials are entirely determined
by the classical differentials.
However, one must be careful with the weights.  We describe this 
issue next.

The first classical differential is a $d_3$ that hits $h_1^4$.
The input to this differential has weight $3$ (since it lies in
the $5$-stem and has filtration $1$),
but $h_1^4$ has weight $4$.  It follows that
$d_3$ hits $\tau h_1^4$ in the motivic calculation, not $h_1^4$
as in the classical calculation.  As a result, instead of finding 
that $h_1^4 = 0$, we have that $h_1^k$ is non-zero for all $k \geq 0$
but $\tau h_1^k = 0$ for $k \geq 4$.

The same situation occurs for the $h_1$-multiples of many of the 
classes in filtration $0$.  Classically, $d_3$-differentials
tell us that these classes are killed by $h_1^3$.  Motivically,
these classes are not killed by $h_1^k$ for any $k$, but they
are killed by $\tau h_1^3$.

In order to simplify our discussion and make our diagrams legible,
from here on we shall ignore these classes in filtration $0$ 
and their $h_1^k$-multiples without further mention.

Classically, the next differential is $d_5(\Delta) = h_2 g$.
Since $\Delta$ has weight $12$ while $h_2 g$ has weight $14$,
the motivic differential is $d_5(\Delta) = \tau^2 h_2 g$.
We observe that $h_2 g$ survives the motivic calculation, but
$\tau^2 h_2 g = 0$.

The entire calculation proceeds similarly.  
Elementary bookkeeping with the weights shows that
a non-zero $d_{2k+1}$-differential hits $\tau^k$ times a generator.
In other words, non-zero $d_{2k+1}$-differentials produce classes
that are killed by $\tau^k$.

The results of the full analysis are shown in the following tables.
These are the $E_{\infty}$-term of a (speculative) spectral sequence
converging to the $2$-complete homotopy groups of the (speculative)
motivic spectrum $\mmf$.

Boxes represent copies of $\Z_2[\tau]$.  Solid dots represent
copies of $\Z/2[\tau]$.  A number $k$ represents $\Z/2[\tau]/\tau^k$.
A horizontal row of symbols at a single location represents
extensions by $2$.  
For example, at location $(3,1)$ there is a copy of $\Z/4[\tau]$,
and at location $(23,5)$ there is a copy of $\Z/4[\tau]/\tau^2$.
More interestingly, at location $(20,4)$
there is a copy of $\Z/8[\tau]$; at location $(40,8)$ there is a copy
of $\Z/8[\tau]/4\tau^2$; and at location $(120,24)$ there is a copy
of $\Z/8[\tau]/(\tau^{11}, 2\tau^6, 4\tau^2)$.

Observe that if we ignore the numbers and just consider the boxes
and solid dots, we obtain the classical picture.  This expresses
the principle that the classical calculation is recovered by
inverting $\tau$.

The blue lines show extensions by $\eta$ and $\nu$ that are detected
in the $E_\infty$-term.  They take generators to generators
in the predictable way.  The blue arrows of slope $1$ indicate infinite
sequences of elements that are connected by $h_1$-multiplications
and are killed by $\tau$.  We have not shown them explicitly
in order to make the diagrams legible.

The red lines show exotic extensions by
$2$, $\eta$, and $\nu$ that are not detected in $E_\infty$.
Beware that these exotic extensions do not take generators to generators.
The first example occurs in the 3-stem, where $4 \nu = \tau \eta^3$.
The next example is $\nu e[25,1] = \tau^2 \epsilon \overline{\kappa}$.
As always, the required power of $\tau$ is easily determined by
the weights.Here we write $e[s,t]$ for the homotopy element that is 
represent by the generator in stem $s$ and filtration $t$.
Note that $\epsilon$ and $\overline{\kappa}$ are names for
$e[8,2]$ and $e[20,4]$.

The motivic calculation shows many more exotic extensions than the
classical calculation.  However, all of the motivic extensions 
are easily implied by the classical ones.  We remark that
the charts in \cite{B} inadvertently failed to indicate the 
exotic extension $\eta e[104,2] = e[105,17]$.

\begin{ex}
Note that $\pi_{0,*} \mmf$ is $\Z_2[\tau]$.  Hence for any $k$,
$\pi_{k,*} \mmf$ is a $\Z_2[\tau]$-module.

Consider $\pi_{120,*} \mmf$.  This module has 6 free generators of weight
$60$.  These arise in filtration $0$; only one is shown on the chart.

In addition, $\pi_{120,*} \mmf$ contains many copies of $\Z/2[\tau]/\tau$;
these arise as $h_1$-multiples and are not shown explicitly on the chart.

Finally, $\pi_{120,*} \mmf$ contains a copy of 
$\Z/8[\tau]/(\tau^{11}, 2\tau^6, 4\tau^2)$, generated by 
$\overline{\kappa}^6$ with weight 72.
\end{ex}

\begin{ex}
The $\Z_2[\tau]$-module $\pi_{170,*} \mmf$ contains 
a copy of 
\[
\Z/8[x,y,z]/ (\tau^2 x, \tau^6 y, \tau^{11} z,
2x = \tau^4 y, 2y = \tau^6 z)
\]
where the generators $x$, $y$, and $z$
have weights 88, 92, and 98 respectively.
\end{ex}

\begin{remark}
In \cite[Proposition 8.4]{B},
the classical relation $\epsilon \kappa = \eta^2 \overline{\kappa}$
is computed by consideration of $\eta^3 \overline{\kappa}$. 
Here $\kappa$ is a name for $e[14,2]$.
However,
the argument is not quite so simple, since this classical product is zero.  
Motivically, this complication does not occur.
We can show that $\epsilon \kappa = \tau \eta^2 \overline{\kappa}$
by considering $\tau \eta^3 \overline{\kappa}$, which is non-zero.
(Note that $\tau^3 \eta^3 \overline{\kappa}$ is zero, but this does
not interfere with the computation.)
\end{remark}

\begin{landscape}

\psset{unit=0.39cm}

\begin{pspicture}(52,12)

\psgrid[unit=2,gridcolor=gridline,subgriddiv=0,gridlabelcolor=white](0,0)(26,6)

\rput(0,-1){\small 0}
\rput(4,-1){\small 4}
\rput(8,-1){\small 8}
\rput(12,-1){\small 12}
\rput(16,-1){\small 16}
\rput(20,-1){\small 20}
\rput(24,-1){\small 24}
\rput(28,-1){\small 28}
\rput(32,-1){\small 32}
\rput(36,-1){\small 36}
\rput(40,-1){\small 40}
\rput(44,-1){\small 44}
\rput(48,-1){\small 48}
\rput(52,-1){\small 52}

\rput(-1,0){\small 0}
\rput(-1,4){\small 4}
\rput(-1,8){\small 8}
\rput(-1,12){\small 12}


\extn(2.7,1)(3,3)
\extn(22.7,5)(23,7)
\extn(42.7,9)(43,11)
\extn(25,1)(28,6)
\extn(27,1)(27,3)
\extn(27,1)(28,6)
\extn(32,2)(35,7)

\extn(50,2)(52,5.3)
\extn(50.7,1)(51,3)
\extn(39,3)(40,8)
\extn(39,3)(42,10)
\extn(45,5)(48,10)
\extn(47,5)(47,7)
\extn(47,5)(48,10)


\mult(0,0)(1,1)
\mult(0,0)(2.7,1)
\mult(1,1)(2,2)
\mult(2,2)(3,3)
\mult(2.7,1)(6,2)
\mult{->}(3,3)(4,4)
\mult(6,2)(9,3)
\mult(8,2)(9,3)

\psframe(-0.4,-0.4)(0.4,0.4)
\rput(1,1){$\bullet$}
\rput(2,2){$\bullet$}
\rput(3,1){$\bullet\bullet$}
\rput(3,3){$\bullet$}
\rput(6,2){$\bullet$}
\rput(8,2){$\bullet$}
\rput(9,3){$\bullet$}


\mult(14,2)(15,3)
\mult(14,2)(17,3)
\mult(17,3)(20.5,4)
\mult(19.5,4)(21,5)
\mult(19.5,4)(22.7,5)
\mult(21,5)(22,6)
\mult(22,6)(23,7)
\mult(22.7,5)(26,6)
\mult{->}(23,7)(24,8)
\mult(26,6)(29,7)
\mult(28,6)(29,7)

\rput(14,2){$\bullet$}
\rput(15,3){$\bullet$}
\rput(17,3){$\bullet$}
\rput(20,4){$\bullet\bullet$$\bullet$}
\rput(21,5){$\bullet$}
\rput(22,6){$\bullet$}
\rput(23,5){\small 22}
\rput(23,7){\small 3}
\rput(26,6){\small 2}
\rput(28,6){$\bullet$}
\rput(29,7){\small 2}


\mult(34,6)(35,7)
\mult(34,6)(37,7)
\mult(37,7)(40.5,8)
\mult(39.5,8)(41,9)
\mult(39.5,8)(42.7,9)
\mult(41,9)(42,10)
\mult(42,10)(43,11)
\mult(42.7,9)(46,10)
\mult{->}(43,11)(44,12)
\mult(46,10)(49,11)
\mult(48,10)(49,11)

\rput(34,6){$\bullet$}
\rput(35,7){$\bullet$}
\rput(37,7){\small 2}
\rput(40,8){$\bullet\bullet$\small 2}
\rput(41,9){$\bullet$}
\rput(42,10){$\bullet$}
\rput(43,9){\small 22}
\rput(43,11){\small 3}
\rput(46,10){\small 2}
\rput(48,10){\small 4}
\rput(49,11){\small 2}


\mult(25,1)(26,2)
\mult(26,2)(27,3)
\mult{->}(27,3)(28,4)
\mult(32,2)(33,3)

\psframe(23.6,-0.4)(24.4,0.4)
\rput(25,1){$\bullet$}
\rput(26,2){$\bullet$}
\rput(27,1){$\bullet$}
\rput(27,3){$\bullet$}
\rput(32,2){$\bullet$}
\rput(33,3){$\bullet$}


\mult(45,5)(46,6)
\mult(46,6)(47,7)
\mult{->}(47,7)(48,8)

\rput(39,3){$\bullet$}
\rput(45,5){$\bullet$}
\rput(46,6){$\bullet$}
\rput(47,5){\small 2}
\rput(47,7){\small 3}
\rput(52,6){$\bullet$}


\psbezier[linecolor=red](48,0)(49,0)(51,2)(51,3)

\mult(50,2)(51,3)
\mult(50.7,1)(52,1.3)
\mult{->}(51,3)(52,4)

\psframe(47.6,-0.4)(48.4,0.4)
\rput(50,2){$\bullet$}
\rput(51,1){$\bullet\bullet$}
\rput(51,3){$\bullet$}

\end{pspicture}


\newpage

\begin{pspicture}(48,0)(100,20)

\psgrid[unit=2,gridcolor=gridline,subgriddiv=0,gridlabelcolor=white](24,0)(50,10)

\rput(48,-1){\small 48}
\rput(52,-1){\small 52}
\rput(56,-1){\small 56}
\rput(60,-1){\small 60}
\rput(64,-1){\small 64}
\rput(68,-1){\small 68}
\rput(72,-1){\small 72}
\rput(76,-1){\small 76}
\rput(80,-1){\small 80}
\rput(84,-1){\small 84}
\rput(88,-1){\small 88}
\rput(92,-1){\small 92}
\rput(96,-1){\small 96}
\rput(100,-1){\small 100}

\rput(47,0){\small 0}
\rput(47,4){\small 4}
\rput(47,8){\small 8}
\rput(47,12){\small 12}
\rput(47,16){\small 16}
\rput(47,20){\small 20}


\extn(62.7,13)(63,15)
\extn(82.7,17)(83,19)
\extn(52,6)(55,11)
\extn(59,7)(60,12)
\extn(59,7)(62,14)
\extn(65,9)(68,14)
\extn(67,9)(67,11)
\extn(67,9)(68,14)
\extn(72,10)(75,15)
\extn(79,11)(80,16)
\extn(79,11)(82,18)
\extn(85,13)(88,18)
\extn(87,13)(87,15)
\extn(87,13)(88,18)
\extn(92,14)(95,19)
\extn(99,15)(100,20)
\extn(99,15)(100,17.3)
\psbezier[linecolor=red](48,0)(49,0)(51,2)(51,3)
\extn(50,2)(53,7)
\extn(50.7,1)(51,3)
\extn(50.7,1)(52,6)
\extn(54,2)(54,10)
\extn(57,3)(60,12)
\extn(65,3)(66,10)
\extn(70,6)(73,11)
\extn(70.7,5)(71,7)
\extn(70.7,5)(72,10)
\extn(74,6)(74,14)
\extn(77,7)(80,16)
\extn(85,7)(86,14)
\psbezier[linecolor=red](88,8)(89,8)(91,10)(91,11)
\extn(90,10)(93,15)
\extn(90.7,9)(91,11)
\extn(90.7,9)(92,14)
\extn(94,10)(94,18)
\extn(97,11)(100,20)
\extn(97,1)(99.5,20)
\extn(98.7,1)(99,3)
\extn(98.7,1)(99.5,20)


\mult(48,10)(49,11)
\mult(48,10.7)(49,11)

\rput(48,10){\small 4}
\rput(49,11){\small 2}


\mult(54,10)(55,11)
\mult(54,10)(57,11)
\mult(57,11)(60.5,12)
\mult(59.5,12)(61,13)
\mult(59.5,12)(62.7,13)
\mult(61,13)(62,14)
\mult(62,14)(63,15)
\mult(62.7,13)(66,14)
\mult{->}(63,15)(64,16)
\mult(66,14)(69,15)
\mult(68,14)(69,15)

\rput(54,10){$\bullet$}
\rput(55,11){\small 4}
\rput(57,11){\small 2}
\rput(60,12){$\bullet\bullet$\small 2}
\rput(61,13){\small 5}
\rput(62,14){\small 5}
\rput(63,13){\small 22}
\rput(63,15){\small 3}
\rput(66,14){\small 2}
\rput(68,14){\small 4}
\rput(69,15){\small 2}


\mult(74,14)(75,15)
\mult(74,14)(77,15)
\mult(77,15)(80.5,16)
\mult(79.5,16)(81,17)
\mult(79.5,16)(82.7,17)
\mult(81,17)(82,18)
\mult(82,18)(83,19)
\mult(82.7,17)(86,18)
\mult{->}(83,19)(84,20)
\mult(86,18)(89,19)
\mult(88,18)(89,19)

\rput(74,14){\small 6}
\rput(75,15){\small 4}
\rput(77,15){\small 2}
\rput(80,16){$\bullet$\small 62}
\rput(81,17){\small 5}
\rput(82,18){\small 5}
\rput(83,17){\small 22}
\rput(83,19){\small 3}
\rput(86,18){\small 2}
\rput(88,18){\small 4}
\rput(89,19){\small 2}


\mult(94,18)(95,19)
\mult(94,18)(97,19)
\mult(97,19)(100.5,20)

\rput(94,18){\small 6}
\rput(95,19){\small 4}
\rput(97,19){\small 2}
\rput(100,20){$\bullet$\small 62}


\mult(52,6)(53,7)

\rput(52,6){$\bullet$}
\rput(53,7){$\bullet$}


\mult(65,9)(66,10)
\mult(66,10)(67,11)
\mult{->}(67,11)(68,12)
\mult(72,10)(73,11)

\rput(59,7){$\bullet$}
\rput(65,9){$\bullet$}
\rput(66,10){$\bullet$}
\rput(67,9){\small 2}
\rput(67,11){\small 3}
\rput(72,10){\small 4}
\rput(73,11){\small 4}


\mult(85,13)(86,14)
\mult(86,14)(87,15)
\mult{->}(87,15)(88,16)
\mult(92,14)(93,15)

\rput(79,11){\small 4}
\rput(85,13){$\bullet$}
\rput(86,14){\small 5}
\rput(87,13){\small 2}
\rput(87,15){\small 3}
\rput(92,14){\small 4}
\rput(93,15){\small 4}


\rput(99,15){\small 4}


\mult(50,2)(51,3)
\mult(50.7,1)(54,2)
\mult{->}(51,3)(52,4)
\mult(54,2)(57,3)

\psframe(47.6,-0.4)(48.4,0.4)
\rput(50,2){$\bullet$}
\rput(51,1){$\bullet\bullet$}
\rput(51,3){$\bullet$}
\rput(54,2){$\bullet$}
\rput(57,3){$\bullet$}


\mult(65,3)(68,4)
\mult(70,6)(71,7)
\mult(70.7,5)(74,6)
\mult{->}(71,7)(72,8)
\mult(74,6)(77,7)

\rput(65,3){$\bullet$}
\rput(68,4){$\bullet$}
\rput(70,6){$\bullet$}
\rput(71,5){\small 22}
\rput(71,7){\small 3}
\rput(74,6){\small 2}
\rput(77,7){\small 2}
\psbezier[linecolor=red](68,4)(69,4)(71,6)(71,7)


\mult(85,7)(88,8)
\mult(90,10)(91,11)
\mult(90.7,9)(94,10)
\mult{->}(91,11)(92,12)
\mult(94,10)(97,11)

\rput(85,7){\small 2}
\rput(88,8){\small 2}
\rput(90,10){$\bullet$}
\rput(91,9){\small 22}
\rput(91,11){\small 3}
\rput(94,10){\small 2}
\rput(97,11){\small 2}


\mult{->}(75,3)(76,4)

\psframe(71.6,-0.4)(72.4,0.4)
\rput(75,3){$\bullet$}


\mult{->}(95,7)(96,8)

\rput(95,7){\small 3}


\mult(96,0)(99.3,1)
\mult(97,1)(98,2)
\mult(98,2)(99,3)
\mult(98.7,1)(100,1.3)
\mult{->}(99,3)(100,4)

\psframe(95.6,-0.4)(96.4,0.4)
\rput(97,1){$\bullet$}
\rput(98,2){$\bullet$}
\rput(99,1){$\bullet\bullet$}
\rput(99,3){$\bullet$}

\end{pspicture}


\begin{pspicture}(96,0)(148,32)

\psgrid[unit=2,gridcolor=gridline,subgriddiv=0,gridlabelcolor=white](48,0)(74,16)

\rput(96,-1){\small 96}
\rput(100,-1){\small 100}
\rput(104,-1){\small 104}
\rput(108,-1){\small 108}
\rput(112,-1){\small 112}
\rput(116,-1){\small 116}
\rput(120,-1){\small 120}
\rput(124,-1){\small 124}
\rput(128,-1){\small 128}
\rput(132,-1){\small 132}
\rput(136,-1){\small 136}
\rput(140,-1){\small 140}
\rput(144,-1){\small 144}
\rput(148,-1){\small 148}

\rput(95,0){\small 0}
\rput(95,4){\small 4}
\rput(95,8){\small 8}
\rput(95,12){\small 12}
\rput(95,16){\small 16}
\rput(95,20){\small 20}
\rput(95,24){\small 24}
\rput(95,28){\small 28}
\rput(95,32){\small 32}


\extn(97,11)(100,20)
\extn(102.7,21)(103,23)
\extn(122.7,25)(123,27)
\extn(142.7,29)(143,31)
\extn(99,15)(102,22)
\extn(99,15)(100,17.3)
\extn(105,17)(108,22)
\extn(107,17)(107,19)
\extn(107,17)(108,22)
\extn(112,18)(115,23)
\extn(119,19)(122,26)
\extn(119,19)(120,21.3)
\extn(125,21)(128,26)
\extn(127,21)(127,23)
\extn(127,21)(128,26)
\extn(132,22)(135,27)
\extn(139,23)(142,30)
\extn(139,23)(140,25.3)
\extn(145,25)(148,30)
\extn(147,25)(147,27)
\extn(147,25)(148,30)
\extn(105,11)(106,18)
\psbezier[linecolor=red](108,12)(109,12)(111,14)(111,15)
\extn(110,14)(113,19)
\extn(110.7,13)(111,15)
\extn(110.7,13)(112,18)
\extn(114,14)(114,22)
\extn(117,15)(120,24)
\extn(125,15)(126,22)
\psbezier[linecolor=red](128,16)(129,16)(131,18)(131,19)
\extn(130,18)(133,23)
\extn(130.7,17)(131,19)
\extn(130.7,17)(132,22)
\extn(134,18)(134,26)
\extn(137,19)(140,28)
\extn(145,19)(146,26)
\extn(97,1)(99.5,20)
\extn(98.7,1)(99,3)
\extn(98.7,1)(99.5,20)
\extn(104,2)(105,17)
\extn(110,2)(110,14)
\extn(117,5)(119.5,24)
\extn(118.7,5)(119,7)
\extn(118.7,5)(119.5,24)
\extn(124,6)(125,21)
\extn(130,6)(130,18)
\extn(137,9)(139.5,28)
\extn(138.7,9)(139,11)
\extn(138.7,9)(139.5,28)
\extn(144,10)(145,25)
\extn(123,1)(123,3)
\extn(123,1)(124,6)
\extn(128,2)(131,7)
\extn(135,3)(135.7,8)
\extn(135,3)(138,10)
\extn(143,5)(143,7)
\extn(143,5)(144,10)
\extn(144,0)(147,3)
\extn(146.7,1)(147,3)
\extn(146.7,1)(148,6)


\mult(96,10.7)(97,11)

\rput(97,11){\small 2}


\mult(96,18.7)(97,19)
\mult(97,19)(100.5,20)
\mult(99.5,20)(101,21)
\mult(99.5,20)(102.7,21)
\mult(101,21)(102,22)
\mult(102,22)(103,23)
\mult(102.7,21)(106,22)
\mult{->}(103,23)(104,24)
\mult(106,22)(109,23)
\mult(108,22)(109,23)

\rput(97,19){\small 2}
\rput(100,20){$\bullet$\small 62}
\rput(101,21){\small 5}
\rput(102,22){\small 5}
\rput(103,21){\small 22}
\rput(103,23){\small 3}
\rput(106,22){\small 2}
\rput(108,22){\small 4}
\rput(109,23){\small 2}


\mult(114,22)(115,23)
\mult(114,22)(117,23)
\mult(116,22.7)(117,23)
\mult(117,23)(120.5,24)
\mult(119.5,24)(121,25)
\mult(119.5,24)(122.7,25)
\mult(121,25)(122,26)
\mult(122,26)(123,27)
\mult(122.7,25)(126,26)
\mult{->}(123,27)(124,28)
\mult(126,26)(129,27)
\mult(128,26)(129,27)

\rput(114,22){\small 6}
\rput(115,23){\small 4}
\rput(117,23){\small 2}
\rput(120,24){\small 11-62}
\rput(121,25){\small 5}
\rput(122,26){\small 5}
\rput(123,25){\small 22}
\rput(123,27){\small 3}
\rput(126,26){\small 2}
\rput(128,26){\small 4}
\rput(129,27){\small 2}


\mult(134,26)(135,27)
\mult(134,26)(137,27)
\mult(136,26.7)(137,27)
\mult(137,27)(140.5,28)
\mult(139.5,28)(141,29)
\mult(139.5,28)(142.7,29)
\mult(141,29)(142,30)
\mult(142,30)(143,31)
\mult(142.7,29)(146,30)
\mult{->}(143,31)(144,32)
\mult(146,30)(148,30.6)

\rput(134,26){\small 6}
\rput(135,27){\small 4}
\rput(137,27){\small 2}
\rput(140,28){\small 11-62}
\rput(141,29){\small 5}
\rput(142,30){\small 5}
\rput(143,29){\small 22}
\rput(143,31){\small 3}
\rput(146,30){\small 2}
\rput(148,30){\small 4}


\mult(105,17)(106,18)
\mult(106,18)(107,19)
\mult{->}(107,19)(108,20)
\mult(112,18)(113,19)

\rput(99,15){\small 4}
\rput(105,17){$\bullet$}
\rput(106,18){\small 5}
\rput(107,17){\small 2}
\rput(107,19){\small 3}
\rput(112,18){\small 4}
\rput(113,19){\small 4}


\mult(125,21)(126,22)
\mult(126,22)(127,23)
\mult{->}(127,23)(128,24)
\mult(132,22)(133,23)

\rput(119,19){\small 4}
\rput(125,21){$\bullet$}
\rput(126,22){\small 5}
\rput(127,21){\small 2}
\rput(127,23){\small 3}
\rput(132,22){\small 4}
\rput(133,23){\small 4}


\mult(145,25)(146,26)
\mult(146,26)(147,27)
\mult{->}(147,27)(148,28)

\rput(139,23){\small 4}
\rput(145,25){\small 11}
\rput(146,26){\small 5}
\rput(147,25){\small 2}
\rput(147,27){\small 3}


\mult(105,11)(108,12)
\mult(110,14)(111,15)
\mult(110.7,13)(114,14)
\mult{->}(111,15)(112,16)
\mult(114,14)(117,15)

\rput(105,11){\small 2}
\rput(108,12){\small 2}
\rput(110,14){$\bullet$}
\rput(111,13){\small 22}
\rput(111,15){\small 3}
\rput(114,14){\small 2}
\rput(117,15){\small 2}


\mult(125,15)(128,16)
\mult(130,18)(131,19)
\mult(130.7,17)(134,18)
\mult{->}(131,19)(132,20)
\mult(134,18)(137,19)

\rput(125,15){\small 2}
\rput(128,16){\small 2}
\rput(130,18){$\bullet$}
\rput(131,17){\small 22}
\rput(131,19){\small 3}
\rput(134,18){\small 2}
\rput(137,19){\small 2}


\mult(145,19)(148,20)

\rput(145,19){\small 2}
\rput(148,20){\small 2}


\mult{->}(115,11)(116,12)

\rput(115,11){\small 3}


\mult{->}(135,15)(136,16)

\rput(135,15){\small 3}


\mult(96,0)(99.3,1)
\mult(97,1)(98,2)
\mult(98,2)(99,3)
\mult(98.7,1)(102,2)
\mult{->}(99,3)(100,4)
\mult(102,2)(105,3)
\mult(104,2)(105,3)

\psframe(95.6,-0.4)(96.4,0.4)
\rput(97,1){$\bullet$}
\rput(98,2){$\bullet$}
\rput(99,1){$\bullet\bullet$}
\rput(99,3){$\bullet$}
\rput(102,2){$\bullet$}
\rput(104,2){$\bullet$}
\rput(105,3){$\bullet$}


\mult(110,2)(111,3)
\mult(110,2)(113,3)
\mult(113,3)(116.3,4)
\mult(115.7,4)(119.3,5)
\mult(117,5)(118,6)
\mult(118,6)(119,7)
\mult(118.7,5)(122,6)
\mult{->}(119,7)(120,8)
\mult(122,6)(125,7)
\mult(124,6)(125,7)

\rput(110,2){$\bullet$}
\rput(111,3){$\bullet$}
\rput(113,3){$\bullet$}
\rput(116,4){$\bullet\bullet$}
\rput(117,5){$\bullet$}
\rput(118,6){$\bullet$}
\rput(119,5){\small 22}
\rput(119,7){\small 3}
\rput(122,6){\small 2}
\rput(124,6){$\bullet$}
\rput(125,7){\small 2}


\mult(130,6)(131,7)
\mult(130,6)(133,7)
\mult(133,7)(136.3,8)
\mult(135.7,8)(139.3,9)
\mult(137,9)(138,10)
\mult(138,10)(139,11)
\mult(138.7,9)(142,10)
\mult{->}(139,11)(140,12)
\mult(142,10)(145,11)
\mult(144,10)(145,11)

\rput(130,6){$\bullet$}
\rput(131,7){$\bullet$}
\rput(133,7){\small 2}
\rput(136,8){$\bullet$\small 2}
\rput(137,9){$\bullet$}
\rput(138,10){$\bullet$}
\rput(139,9){\small 22}
\rput(139,11){\small 3}
\rput(142,10){\small 2}
\rput(144,10){\small 4}
\rput(145,11){\small 2}


\mult(122,2)(123,3)
\mult{->}(123,3)(124,4)
\mult(128,2)(129,3)

\psframe(119.6,-0.4)(120.4,0.4)
\rput(122,2){$\bullet$}
\rput(123,1){$\bullet$}
\rput(123,3){$\bullet$}
\rput(128,2){$\bullet$}
\rput(129,3){$\bullet$}


\mult(142,6)(143,7)
\mult{->}(143,7)(144,8)

\rput(135,3){$\bullet$}
\rput(142,6){$\bullet$}
\rput(143,5){\small 2}
\rput(143,7){\small 3}
\rput(148,6){$\bullet$}


\mult(146.7,1)(148,1.3)
\mult{->}(147,3)(148,4)

\psframe(143.6,-0.4)(144.4,0.4)
\rput(147,1){$\bullet\bullet$}
\rput(147,3){$\bullet$}

\end{pspicture}


\psset{unit=0.36cm}

\begin{pspicture}(144,0)(196,40)

\psgrid[unit=2,gridcolor=gridline,subgriddiv=0,gridlabelcolor=white](72,0)(98,20)

\rput(144,-1){\small 144}
\rput(148,-1){\small 148}
\rput(152,-1){\small 152}
\rput(156,-1){\small 156}
\rput(160,-1){\small 160}
\rput(164,-1){\small 164}
\rput(168,-1){\small 168}
\rput(172,-1){\small 172}
\rput(176,-1){\small 176}
\rput(180,-1){\small 180}
\rput(184,-1){\small 184}
\rput(188,-1){\small 188}
\rput(192,-1){\small 192}
\rput(196,-1){\small 196}

\rput(143,0){\small 0}
\rput(143,4){\small 4}
\rput(143,8){\small 8}
\rput(143,12){\small 12}
\rput(143,16){\small 16}
\rput(143,20){\small 20}
\rput(143,24){\small 24}
\rput(143,28){\small 28}
\rput(143,32){\small 32}
\rput(143,36){\small 36}
\rput(143,40){\small 40}


\extn(162.7,33)(163,35)
\extn(182.7,37)(183,39)
\extn(145,25)(148,30)
\extn(147,25)(147,27)
\extn(147,25)(148,30)
\extn(152,26)(155,31)
\extn(179,31)(182,38)
\extn(179,31)(180,33.3)
\extn(185,33)(188,38)
\extn(187,33)(187,35)
\extn(187,33)(188,38)
\extn(192,34)(195,39)
\extn(145,19)(146,26)
\psbezier[linecolor=red](148,20)(149,20)(151,22)(151,23)
\extn(150,22)(153,27)
\extn(150.7,21)(151,23)
\extn(150.7,21)(152,26)
\extn(154,22)(154,30)
\extn(157,23)(160,32)
\extn(159,27)(162,34)
\extn(159,27)(160,29.3)
\extn(165,29)(168,34)
\extn(167,29)(167,31)
\extn(167,29)(168,34)
\extn(172,30)(175,35)
\extn(165,23)(166,30)
\psbezier[linecolor=red](168,24)(169,24)(171,26)(171,27)
\extn(170,26)(173,31)
\extn(170.7,25)(171,27)
\extn(170.7,25)(172,30)
\extn(174,26)(174,34)
\extn(177,27)(180,36)
\extn(185,27)(186,34)
\psbezier[linecolor=red](188,28)(189,28)(191,30)(191,31)
\extn(190,30)(193,35)
\extn(190.7,29)(191,31)
\extn(190.7,29)(192,34)
\extn(194,30)(194,38)
\extn(150,10)(150,22)
\extn(157,13)(159.5,32)
\extn(158.7,13)(159,15)
\extn(158.7,13)(159.5,32)
\extn(164,14)(165,29)
\extn(170,14)(170,26)
\extn(177,17)(179.5,36)
\extn(178.7,17)(179,19)
\extn(178.7,17)(179.5,36)
\extn(184,18)(185,33)
\extn(190,18)(190,30)
\extn(168,10)(171,15)
\extn(148,6)(151,11)
\extn(155,7)(155.7,12)
\extn(155,7)(158,14)
\extn(163,9)(163,11)
\extn(163,9)(164,14)
\extn(175,11)(175.7,16)
\extn(175,11)(178,18)
\extn(183,13)(183,15)
\extn(183,13)(184,18)
\extn(188,14)(191,19)
\extn(195,15)(195.7,20)
\extn(195,15)(196,17.3)
\extn(144,0)(147,3)
\extn(146.7,1)(147,3)
\extn(146.7,1)(148,6)
\extn(150,2)(150,10)
\extn(153,3)(155.7,12)
\extn(161,3)(162,10)
\extn(164,4)(167,7)
\extn(166.7,5)(167,7)
\extn(166.7,5)(168,10)
\extn(170,6)(170,14)
\extn(173,7)(175.7,16)
\extn(181,7)(182,14)
\extn(184,8)(187,11)
\extn(186.7,9)(187,11)
\extn(186.7,9)(188,14)
\extn(190,10)(190,18)
\extn(193,11)(195.7,20)


\mult(144,10.6)(145,11)
\mult(144,10)(145,11)

\rput(144,10){\small 4}
\rput(145,11){\small 2}


\mult(144,29.3)(146,30)
\mult(146,30)(149,31)
\mult(148,30)(149,31)

\rput(146,30){\small 2}
\rput(148,30){\small 4}
\rput(149,31){\small 2}


\mult(154,30)(155,31)
\mult(154,30)(157,31)
\mult(156,30.7)(157,31)
\mult(157,31)(160.5,32)
\mult(159.5,32)(161,33)
\mult(159.5,32)(162.7,33)
\mult(161,33)(162,34)
\mult(162,34)(163,35)
\mult(162.7,33)(166,34)
\mult{->}(163,35)(164,36)
\mult(166,34)(169,35)
\mult(168,34)(169,35)

\rput(154,30){\small 6}
\rput(155,31){\small 4}
\rput(157,31){\small 2}
\rput(160,32){\small 11-62}
\rput(161,33){\small 5}
\rput(162,34){\small 5}
\rput(163,33){\small 22}
\rput(163,35){\small 3}
\rput(166,34){\small 2}
\rput(168,34){\small 4}
\rput(169,35){\small 2}


\mult(174,34)(175,35)
\mult(174,34)(177,35)
\mult(176,34.7)(177,35)
\mult(177,35)(180.5,36)
\mult(179.5,36)(181,37)
\mult(179.5,36)(182.7,37)
\mult(181,37)(182,38)
\mult(182,38)(183,39)
\mult(182.7,37)(186,38)
\mult{->}(183,39)(184,40)
\mult(186,38)(189,39)
\mult(188,38)(189,39)

\rput(174,34){\small 6}
\rput(175,35){\small 4}
\rput(177,35){\small 2}
\rput(180,36){\small 11-62}
\rput(181,37){\small 5}
\rput(182,38){\small 5}
\rput(183,37){\small 22}
\rput(183,39){\small 3}
\rput(186,38){\small 2}
\rput(188,38){\small 4}
\rput(189,39){\small 2}


\mult(194,38)(195,39)
\mult(194,38)(196,38.6)

\rput(194,38){\small 6}
\rput(195,39){\small 4}


\mult(145,25)(146,26)
\mult(146,26)(147,27)
\mult{->}(147,27)(148,28)
\mult(152,26)(153,27)

\rput(145,25){\small 11}
\rput(146,26){\small 5}
\rput(147,25){\small 2}
\rput(147,27){\small 3}
\rput(152,26){\small 4}
\rput(153,27){\small 4}


\mult{->}(167,31)(168,32)
\mult(172,30)(173,31)
\mult(165,29)(166,30)
\mult(166,30)(167,31)

\rput(159,27){\small 4}
\rput(165,29){\small 11}
\rput(166,30){\small 5}
\rput(167,29){\small 2}
\rput(167,31){\small 3}
\rput(172,30){\small 4}
\rput(173,31){\small 4}


\mult(185,33)(186,34)
\mult(186,34)(187,35)
\mult{->}(187,35)(188,36)
\mult(192,34)(193,35)

\rput(179,31){\small 4}
\rput(185,33){\small 11}
\rput(186,34){\small 5}
\rput(187,33){\small 2}
\rput(187,35){\small 3}
\rput(192,34){\small 4}
\rput(193,35){\small 4}


\mult(145,19)(148,20)
\mult(150,22)(151,23)
\mult(150.7,21)(154,22)
\mult{->}(151,23)(152,24)
\mult(154,22)(157,23)

\rput(145,19){\small 2}
\rput(148,20){\small 2}
\rput(150,22){$\bullet$}
\rput(151,21){\small 22}
\rput(151,23){\small 3}
\rput(154,22){\small 2}
\rput(157,23){\small 2}


\mult(165,23)(168,24)
\mult(170,26)(171,27)
\mult(170.7,25)(174,26)
\mult{->}(171,27)(172,28)
\mult(174,26)(177,27)

\rput(165,23){\small 2}
\rput(168,24){\small 2}
\rput(170,26){\small 11}
\rput(171,25){\small 22}
\rput(171,27){\small 3}
\rput(174,26){\small 2}
\rput(177,27){\small 2}


\mult(185,27)(188,28)
\mult(190,30)(191,31)
\mult(190.7,29)(194,30)
\mult{->}(191,31)(192,32)
\mult(194,30)(196,30.6)

\rput(185,27){\small 2}
\rput(188,28){\small 2}
\rput(190,30){\small 11}
\rput(191,29){\small 22}
\rput(191,31){\small 3}
\rput(194,30){\small 2}


\mult{->}(155,19)(156,20)

\rput(155,19){\small 3}


\mult{->}(175,23)(176,24)

\rput(175,23){\small 3}


\mult{->}(195,27)(196,28)

\rput(195,27){\small 3}


\mult(150,10)(151,11)
\mult(150,10)(153,11)
\mult(153,11)(156.3,12)
\mult(155.7,12)(159.3,13)
\mult(157,13)(158,14)
\mult(158,14)(159,15)
\mult(158.7,13)(162,14)
\mult{->}(159,15)(160,16)
\mult(162,14)(165,15)
\mult(164,14)(165,15)

\rput(150,10){$\bullet$}
\rput(151,11){\small 4}
\rput(153,11){\small 2}
\rput(156,12){$\bullet$\small 2}
\rput(157,13){\small 5}
\rput(158,14){\small 5}
\rput(159,13){\small 22}
\rput(159,15){\small 3}
\rput(162,14){\small 2}
\rput(164,14){\small 4}
\rput(165,15){\small 2}


\mult(170,14)(171,15)
\mult(170,14)(173,15)
\mult(173,15)(176.3,16)
\mult(175.7,16)(179.3,17)
\mult(177,17)(178,18)
\mult(178,18)(179,19)
\mult{->}(179,19)(180,20)
\mult(182,18)(185,19)
\mult(184,18)(185,19)
\mult(178.7,17)(182,18)

\rput(170,14){\small 6}
\rput(171,15){\small 4}
\rput(173,15){\small 2}
\rput(176,16){\small 62}
\rput(177,17){\small 5}
\rput(178,18){\small 5}
\rput(179,17){\small 22}
\rput(179,19){\small 3}
\rput(182,18){\small 2}
\rput(184,18){\small 4}
\rput(185,19){\small 2}


\mult(190,18)(191,19)
\mult(190,18)(193,19)
\mult(193,19)(196.3,20)

\rput(196,20){\small 62}
\rput(190,18){\small 6}
\rput(191,19){\small 4}
\rput(193,19){\small 2}


\mult(148,6)(149,7)

\rput(148,6){$\bullet$}
\rput(149,7){$\bullet$}


\mult(162,10)(163,11)
\mult{->}(163,11)(164,12)
\mult(168,10)(169,11)

\rput(155,7){$\bullet$}
\rput(162,10){$\bullet$}
\rput(163,9){\small 2}
\rput(163,11){\small 3}
\rput(168,10){\small 4}
\rput(169,11){\small 4}


\mult(182,14)(183,15)
\mult{->}(183,15)(184,16)
\mult(188,14)(189,15)

\rput(175,11){\small 4}
\rput(182,14){\small 5}
\rput(183,13){\small 2}
\rput(183,15){\small 3}
\rput(188,14){\small 4}
\rput(189,15){\small 4}


\rput(195,15){\small 4}


\mult(146.7,1)(150,2)
\mult{->}(147,3)(148,4)
\mult(150,2)(153,3)

\psframe(143.6,-0.4)(144.4,0.4)
\rput(147,1){$\bullet\bullet$}
\rput(147,3){$\bullet$}
\rput(150,2){$\bullet$}
\rput(153,3){$\bullet$}


\mult(161,3)(164,4)
\mult(166.7,5)(170,6)
\mult{->}(167,7)(168,8)
\mult(170,6)(173,7)

\rput(161,3){$\bullet$}
\rput(164,4){$\bullet$}
\rput(167,5){\small 22}
\rput(167,7){\small 3}
\rput(170,6){\small 2}
\rput(173,7){\small 2}


\mult(181,7)(184,8)
\mult(186.7,9)(190,10)
\mult{->}(187,11)(188,12)
\mult(190,10)(193,11)

\rput(181,7){\small 2}
\rput(184,8){\small 2}
\rput(187,9){\small 22}
\rput(187,11){\small 3}
\rput(190,10){\small 2}
\rput(193,11){\small 2}


\psframe(167.6,-0.4)(168.4,0.4)

\end{pspicture}

\end{landscape}


\bibliographystyle{amsalpha}

\end{document}